\input ./style/arxiv-general.cfg
\documentclass[aap,MSNbibl,seceqn,dvips]{arximspdf}
\makeatletter
   \@ifpackageloaded{graphicx}{}{\usepackage{graphicx}}
\makeatother

\doi{10.1214/15-AAP1117}
\volume{26}
\issue{3}
\pubyear{2016}
\firstpage{1297}
\lastpage{1328}
\docsubty{FLA}

\makeatletter
\newcommand{\rrvert}{\vert}
\newcommand{\llvert}{\vert}
\newcommand{\eqref}[1]{(\ref{#1})}
\newproclaim{definition}{Definition}[section]
\newtheorem{theorem}{Theorem}[section]
\newtheorem{proposition}{Proposition}[section]
\newtheorem{lemma}{Lemma}[section]
\newproclaim{remark}{Remark}[section]
\newcommand{\var}{\operatorname{Var}}
\newcommand{\real}{\operatorname{Re}}
\makeatother

\begin{document}
\begin{frontmatter}

\title{Social contact processes and the partner model}
\runtitle{The partner model}

\begin{aug}
\author[A]{\fnms{Eric}~\snm{Foxall}\corref{}\thanksref{T1}\ead
[label=e1]{efoxall@uvic.ca}},
\author[A]{\fnms{Roderick}~\snm{Edwards}\thanksref{T2}\ead[label=e2]{edwards@uvic.ca}}~\and
\author[A]{\fnms{P.}~\snm{van~den~Driessche}\thanksref{T2}\ead
[label=e3]{vandendr@uvic.ca}}
\runauthor{E. Foxall, R. Edwards and P. van den Driessche}
\affiliation{University of Victoria}
\address[A]{Department of Mathematics and Statistics\\
University of Victoria \\
P.O. BOX 1700 STN CSC\\
Victoria, British Columbia V8W 2Y2\\
Canada\\
\printead{e1}\\
\phantom{E-mail:\ }\printead*{e2}\\
\phantom{E-mail:\ }\printead*{e3}}
\end{aug}
\thankstext{T1}{Supported in part by an NSERC PGSD2 Graduate Scholarship.}
\thankstext{T2}{Supported in part by NSERC Discovery Grants.}

%
\received{\smonth{12} \syear{2014}}
%
\revised{\smonth{4} \syear{2015}}

%
\begin{abstract}
We consider a stochastic model of infection spread on the complete
graph on $N$ vertices incorporating dynamic partnerships, which we
assume to be monogamous.
This can be seen as a variation on the contact process in which some
form of edge dynamics determines the set of contacts at each moment in time.
We identify a basic reproduction number $R_0$ with the property that if
$R_0<1$ the infection dies out by time $O(\log N)$, while if $R_0>1$
the infection survives
for an amount of time $e^{\gamma N}$ for some $\gamma>0$ and hovers
around a uniquely determined metastable proportion of infectious individuals.
The proof in both cases relies on comparison to a set of mean-field
equations when the infection is widespread, and to a branching process
when the infection is sparse.
\end{abstract}

%
\end{frontmatter}

\section{Introduction}
The contact process is a well-studied model of the spread of an
infection, in which an undirected graph $G=(V,E)$ determines a
collection of sites $V$ and edges $E$ which we can think of as
individuals and as links between individuals along which the infection
can be transmitted. Each site is either healthy or infectious;
infectious sites recover at a certain fixed rate, which is usually
normalized to $1$, and transmit the infection to each of their
neighbours at rate $\lambda$.

The contact process has been studied in a variety of different
settings, including lattices \cite{speed,crit,ips,sis} (to cite just a few), infinite trees \cite{trees}, power law
graphs \cite{plg,mvy} and complete graphs \cite{comp}. In each
case, there is a critical value $\lambda_c$ below which the infection
quickly vanishes from the graph, and above which the infection has a
positive probability of surviving either for all time (if the graph is
infinite), or for an amount of time that grows quickly (either
exponentially or at least faster than polynomially) with the size of
the graph; in the power law case $\lambda_c=0$ so long-time survival is
possible whenever $\lambda>0$.

In a social context, $G$ might describe a contact network in which an
edge connects sites $x$ and $y$ if and only if the corresponding
individuals have sufficiently frequent interactions that infection can
be spread from one to the other. In the contact process, the contact
network is fixed, that is, a given pair of individuals is either
connected or not connected for all time. However, we can easily imagine
a scenario in which connections form and break up dynamically, which we
can model by having edges open and close according to certain rules;
here, we use the convention of percolation theory, in which ``open''
means there is a connection across the edge; note this is the opposite
of the convention for electric circuits. In this case, the edges $E$
represent \emph{possible} connections and we have a process
$E_t\subseteq E$ that describes the set of open edges as a function of
time. This type of process we will call a \emph{social contact
process}, since it involves some form of social dynamics.

In the simplest case, edges open and close independently at some fixed
rates $r_+$ and $r_-$. In this case, the distribution of open edges at
a given time converges to the product measure on $\{0,1\}^E$ with
density $r_+/(r_-+r_+)$. Estimates on the survival region can then be
obtained using the results of \cite{broman} and following the pattern
of~\cite{remenik}. On the other hand, edge dynamics could depend on the
state of the infection; for example, site $x$ might be less likely to
connect with site $y$, if $y$ is infected. If we then relax the
tendency to avoid infected sites, then for a given value of $\lambda$,
we might ask at what point does the infection start to spread, if it
does.

Here, we consider edges opening and closing independently as described
above but with the added restriction of \emph{monogamy}, that is, if
two sites are connected (i.e., linked by an edge) then so long as they
remain connected, they cannot connect to other sites. In this model, we
think of connected pairs as partners, so we call it the \emph{partner
model}. For simplicity, we study the model on the sequence of complete
graphs $K_N$ on $N$ vertices, where $N$ will tend to $\infty$; this is
a reasonable model for, say, the spread of a sexually transmitted
infection through a population of monogamous homosexual individuals in
a big city. We rescale the partner formation rate per edge to $r_+/N$
to ensure that a given individual in a pool of entirely singles finds a
partner at total rate approximately $r_+$. For future reference, we use
interchangeably both the words healthy and susceptible, and the words
unpartnered and single, to describe respectively an individual that is
not infectious, or an individual that does not have a partner. Even in
this simple model, as described below, there is a phase transition
between extinction and spread of the infection.

\section{Statement of main results}\label{secmain}
In order to analyze the partner model, we should first ensure that it
is well defined, so following \cite{gc} we give a graphical
construction which makes it easy to visualize its evolution in time and
space. We write the model as $(V_t,E_t)$ where $V_t \subseteq V$ is the
set of infectious sites at time $t$ and $E_t \subseteq E$ is the set of
open edges at time $t$. In general, we assume $\min(r_+,r_-,\lambda)
>0$ since if any of the parameters is equal to zero the dynamics are
trivial.

The complete graph $K_N=(V,E)$ has sites $V=\{1,\ldots,N\}$ and edges $E=\{
\{x,y\}:x,y \in\{\{1,\ldots,N\},x \neq y\}$. On the spacetime set
$K_N\times[0,\infty)$, place independent Poisson point processes
(p.p.p.s) along the fibers $\{\cdot\}\times[0,\infty)$ as follows:
\begin{itemize}
\item for recovery, at each site with intensity $1$ and label $\times$,
\item for transmission, along each edge $xy \in E$ with intensity
$\lambda$ and label $\leftrightarrow$,
\item for partnership formation, along each edge with intensity $r_+/N$
and label $\uparrow$,
\item for partnership breakup, along each edge with intensity $r_-$ and
label $\downarrow$.
\end{itemize}
These define the probability space $\Omega$, whose realizations $\omega
\in\Omega$ consist of collections of labelled points on $K_N\times
[0,\infty)$. Since the graph is finite, the total intensity of p.p.p.s
is finite, thus with probability 1 events are well ordered in time.
Fixing an admissible initial configuration $(V_0,E_0)$, that is, such
that no two edges $xy$ and $yz$ are both open, we determine $(V_t,E_t)$
as follows. For a well-ordered realization with event times
$t_1<t_2<t_3<\cdots,$ suppose $(V_{t_i},E_{t_i})$ is known. If the event
at time $t_{i+1}$ is:
\begin{itemize}
\item an $\times$ at site $x$ and $x \in V_{t_i}$ then $V_{t_{i+1}} =
V_{t_i}\setminus\{x\}$,
\item a $\leftrightarrow$ along edge $xy$, $xy \in E_{t_i}$, $x \in
V_{t_i}$ and $y \notin V_{t_i}$ then $V_{t_{i+1}} = V_{t_i}\cup\{y\}$,
\item a $\uparrow$ along edge $xy$ and $xz,zy \notin E_{t_i}$ for all
$z$ then $E_{t_{i+1}} = E_{t_i} \cup\{xy\}$,
\item a $\downarrow$ along edge $xy$ and $xy \in E_{t_i}$ then
$E_{t_{i+1}} = E_{t_i} \setminus\{xy\}$.
\end{itemize}
Otherwise the configuration is unchanged. This gives $(V_t,E_t)$ at
times $t_0:=0,t_1,t_2,\ldots;$ for $t \in(t_i,t_{i+1})$ set $V_t =
V_{t_i}$ and $E_t = E_{t_i}$.

For the partner model, we are mostly concerned not with the exact
values of $V_t$ and $E_t$ but with the total number of susceptible and
infectious singles $S_t$ and $I_t$ and the total number of partnered
pairs $\mathit{SS}_t,\mathit{SI}_t,\mathit{II}_t$ of the three possible types; as shown in Section~\ref{secmf}, for each $N$, $(S_t,I_t,\mathit{SS}_t,
\mathit{SI}_t,\mathit{II}_t)$ is a continuous
time Markov chain. In general, it will be more convenient to work with
the rescaled quantities $s_t=S_t/N$, $i_t=I_t/N$, $ss_t=\mathit{SS}_t/N$,
$si_t=\mathit{SI}_t/N$ and $ii_t=\mathit{II}_t/N$.

\begin{figure}

\includegraphics{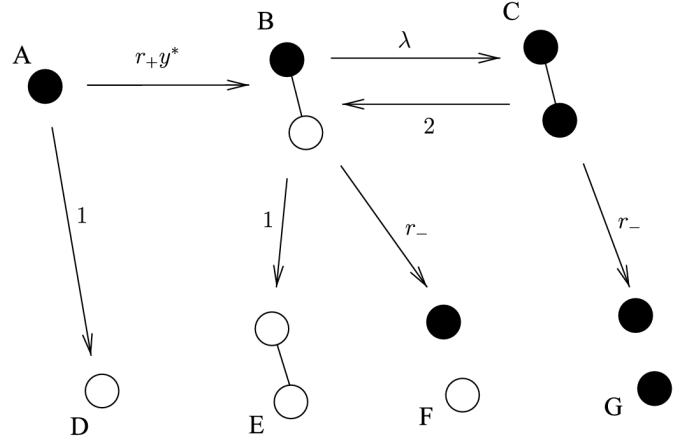}

\caption{Markov chain used to compute $R_0$, with transition rates
indicated; infectious sites are shaded.}
\label{figri}
\end{figure}

Starting from any configuration, as shown in Section~\ref{secstoch},
after a short time the proportion of singles $y_t:=s_t+i_t$ approaches
and remains close to a certain fixed value $y^* \in(0,1)$. The
computation of $y^*$ is given in Section~\ref{secedge}: setting $\alpha
= r_+/r_-$, we find that
%
\begin{equation}
\label{ystar} y^*=1/(2\alpha)[-1 + \sqrt{1 + 4\alpha}].
\end{equation}

To determine the conditions under which the infection can spread, we
use a heuristic argument. Once we know the correct values, we can then
worry about proving they are correct. Suppose we start with $V_0 = \{x\}
$ for some $x \in V$ with $x$ single and $y_0 \approx y^*$, and keep
track of $x$ until the first moment when $x$ either:
\begin{itemize}
\item recovers without finding a partner, or
\item if it finds a partner before recovering, breaks up from that partnership.
\end{itemize}
This leads to the continuous time Markov chain shown in Figure~\ref{figri}. Each of
$A,B,\ldots,G$ represents a state for the chain, and
arrows show possible transitions, with the arrow labelled by the
transition rate. Shaded circles represent infectious individuals and
unshaded circles, healthy individuals. A pair of circles connected by a
line represents a partnered pair. Starting from $A$, a single
infectious site either recovers (goes to $D$) at rate $1$, or finds a
healthy partner at rate $r_+y^*$. Infection takes place at rate $\lambda$. If only one individual in a partnership is infectious (state $B$),
then it recovers at rate 1 (state $E$), and we do not need to worry
about them any more, since neither is infectious. If both are
infectious (state $C$), recovery of one or the other occurs at rate
$2$. While in a partnership, breakup occurs at rate $r_-$.

Define the \emph{basic reproduction number}
%
\begin{equation}
\label{eqri} R_0 = \mathbb{P}(A\rightarrow F) + 2\mathbb{P}(A
\rightarrow G)
\end{equation}
which is the expected number of infectious singles upon absorption of
the above Markov chain, starting from state $A$. As intuition suggests,
and Theorem~\ref{thm1} confirms, the infection can spread if $R_0>1$,
and cannot spread if $R_0 \leq1$.

If the dynamics is in equilibrium, that is, $(s_t,i_t,ss_t,si_t,ii_t)$
hovers around a~fixed value $(s^*,i^*,ss^*,si^*,ii^*)$, then in
particular the proportion of infectious singles is roughly constant. To
compute this proportion, we again use a heuristic argument. Three
events affect infectious singles:
\begin{itemize}
\item $I\rightarrow S$, which occurs at rate $I_t = i_tN$,
\item $I+I\rightarrow \mathit{II}$, which occurs\vspace*{1pt} at rate $(r_+/N){I_t \choose 2}
\approx r_+(i_t^2/2)N$, and
\item $S+I\rightarrow \mathit{SI}$, which occurs at rate $(r_+/N)I_tS_t = r_+i_ts_tN$.
\end{itemize}
If a partnership is formed, then using these rates and Figure~\ref{figri}, we can compute the expected number of infectious singles upon
breakup. Fixing $i_t=i$ for some $i \in[0,y^*]$ and $s_t+i_t=y^*$,
define the normalizing constant $z = 1 + r_+i/2 + r_+(y^*-i) = 1 +
r_+(y^*-i/2)$ and the probabilities $p_S = 1/z$, $p_{\mathit{II}} = r_+i/(2z)$
and $p_{\mathit{SI}} = r_+(y^*-i)/z$ and let
%
\begin{equation}
\label{eqdi} \Delta(i) = p_S\Delta_S +
p_{\mathit{II}}\Delta_{\mathit{II}} + p_{\mathit{SI}}\Delta_{\mathit{SI}},
\end{equation}
where $\Delta_S=-1$, $\Delta_{\mathit{II}} = -2 + \mathbb{P}(C\rightarrow F) +
2\mathbb{P}(C\rightarrow G)$ and $\Delta_{\mathit{SI}} = -1 + \mathbb
{P}(B\rightarrow F) + 2\mathbb{P}(B\rightarrow G)$. The function $\Delta
(i)$ tracks the expected change in the number of infectious singles,
per event affecting one or more infectious singles. Thus, for an
equilibrium solution we should have $\Delta(i^*)=0$. As shown in Lemma~\ref{deltalemma}, to have a solution with $i^*>0$, we need $R_0>1$.

As shown in Lemma~\ref{r0up}, for fixed $r_+,r_-$, $R_0$ is continuous
and increasing in $\lambda$. Defining
%
\begin{equation}
\label{eqlc} \lambda_c = \sup\{\lambda\geq0:R_0 \leq1
\}
\end{equation}
with $\sup\mathbb{R}_+ :=\infty$, it follows that if $\lambda_c=\infty
$ then $R_0<1$ for all $\lambda$, and if $\lambda_c<\infty$ then
$R_0<1$ if $\lambda<\lambda_c$, $R_0=1$ if $\lambda=\lambda_c$ and
$R_0>1$ if $\lambda>\lambda_c$. In models exhibiting a phase
transition, one often seeks a \emph{critical exponent} $\gamma$ such
that for an observable $F(\lambda)$ it holds that $F(\lambda) \sim
C(\lambda-\lambda_c)^{\gamma}$. As we see in the statement of the
upcoming Theorem~\ref{thm2}, here the critical exponent for $i^*$ is
equal to $1$.

The following two theorems are the main results of this paper. The
first result tells us where and when we should expect a phase
transition to occur. In particular, it gives a formula for $\lambda_c$
and describes the behaviour of $i^*$ near $\lambda_c$.

\begin{theorem}\label{thm2}
Let $y^*,R_0,\Delta(i)$ and $\lambda_c$ be as in \eqref{ystar}, \eqref
{eqri}, \eqref{eqdi} and \eqref{eqlc} and let $r_+,r_-$ be fixed. Then
$\lambda_c<\infty\Leftrightarrow r_+y^*>1 \Leftrightarrow r_+ >
1+1/r_-$ and in this case
\begin{eqnarray*}
&& \lambda_c = \frac{2}{r_-}\frac{2}{(r_+y^*-1)} +
\frac{2}{r_-} + \frac
{4}{r_+y^*-1} + 1 + \frac{r_-}{r_+y^*-1}.
\end{eqnarray*}
If $R_0=R_0(\lambda)>1$, there is a unique solution $i^*(\lambda)\in
(0,y^*)$ to the equation $\Delta(i^*)=0$ and $i^*(\lambda) \sim
C(\lambda-\lambda_c)$ as $\lambda\downarrow\lambda_c$, for some
constant $C>0$.
\end{theorem}

The second result shows that our heuristics are correct. More
precisely, $R_0>1$ is a necessary and sufficient condition for spread
and long-time survival of the infection. Moreover, when $R_0>1$ there
is a unique and globally stable endemic equilibrium with $i^*>0$ given
by $\Delta(i^*)=0$.

\begin{theorem}\label{thm1}
Fix $\lambda,r_+,r_-$ and let $y^*, R_0$ and $\Delta(i)$ be as defined
in \eqref{ystar}, \eqref{eqri} and \eqref{eqdi}.
\begin{itemize}
\item If $R_0\leq1$, for each $\varepsilon>0$ there are constants
$C,T,\gamma>0$ so that, from any initial configuration, with
probability $\geq1-Ce^{-\gamma N}$, $|V_T|\leq\varepsilon N$.
\item If $R_0<1$ there are constants $C,T,\gamma>0$ so that, from any
initial configuration, with probability tending to $1$ as $N\rightarrow
\infty$ all sites are healthy by time $T+C \log N$.
\item If $R_0>1$, there is a unique vector $(s^*,i^*,ss^*,si^*,ii^*)$,
satisfying $i^*>0$, $s^*+i^*=y^*$ and $\Delta(i^*)=0$, such that
\begin{itemize}
\item for\vspace*{1pt} each $\varepsilon>0$, there are constants $C,T,\gamma>0$ so
that, from any initial configuration with $|V_0| \geq\varepsilon N$, with\vspace*{1pt}
probability $\geq1-Ce^{-\gamma N}$,
$|(s_t,i_t,ss_t,\break si_t,ii_t)-(s^*,i^*,ss^*,si^*,ii^*)|\leq\varepsilon$ for
$T\leq t \leq e^{\gamma N}$, and
\item there are constants $\delta,p,C,T>0$ so that, from any initial
configuration with $|V_0|>0$, with probability $\geq p$, $|V_{T+C\log
N}|\geq\delta N$.
\end{itemize}
\end{itemize}
\end{theorem}

To obtain the value of the endemic equilibrium and the behaviour when
$|V_0| \geq\varepsilon N$, which we call the \emph{macroscopic} regime,
we use the \emph{mean-field equations} (MFE) introduced in Section~\ref{secmf}, which are a set of differential equations that give a good
approximation to the evolution of $(s_t,i_t,ss_t,si_t,ii_t)$ when $N$
is large. To describe the behaviour when $1\leq|V_0|\leq\varepsilon N$
for small $\varepsilon>0$, which we call the \emph{microscopic} regime, we
use comparison to a branching process; if $R_0<1$ we bound above and if
$R_0>1$ we bound below.

The paper is laid out as follows. Sections~\ref{secedge} and \ref{secsurv} contain the heuristic calculations that allow us to determine
$y^*,R_0,\lambda_c,\Delta(i)$ and prove Theorem~\ref{thm2}. In Section~\ref{secedge}, we give an informal description of the edge dynamics and
compute~$y^*$. In Section~\ref{secsurv}, we analyze $R_0,\lambda
_c,\Delta(i)$ and prove Theorem~\ref{thm2}, in two parts: Propositions~\ref{thm2.1} and~\ref{thm2.2}. In Section~\ref{secmf}, we
introduce the mean-field equations and characterize their dynamics. In
Sections~\ref{secstoch}, \ref{secmacro} and \ref{secmicro}, we consider
the stochastic process and prove Theorem~\ref{thm1}. In Section~\ref{secstoch}, we develop the tools needed to relate the stochastic model
to the mean-field equations. In Section~\ref{secmacro}, we prove the
macroscopic part of Theorem~\ref{thm1}, and in Section~\ref{secmicro}
we prove the microscopic part.

\section{Proportion of singles}\label{secedge}
Starting from the total number of singles $Y_t=S_t+I_t$, the
transitions are:
\begin{itemize}
\item $Y \rightarrow Y-2$ at rate $(r_+/N)Y(Y-1)/2$,
\item $Y \rightarrow Y+2$ at rate $(N-Y)r_-/2$,
\end{itemize}
which for $y_t := Y_t/N$ gives:
\begin{itemize}
\item $y\rightarrow y-2/N$ at rate $[r_+y(y-1/N)/2]N = (r_+y^2/2)N -
r_+y/2$,
\item $y\rightarrow y+2/N$ at rate $[(1-y)r_-/2]N$.
\end{itemize}
Combining these transitions gives
\begin{eqnarray*}
&& \frac{d}{dt}\mathbb{E}(y_t \vert y_t=y) =
-r_+y^2 + r_-(1-y) + \frac{r_+y}{N}.
\end{eqnarray*}
In Lemma~\ref{ydyn}, we make a rigorous statement about the behaviour
of $y_t$. For now, though, some heuristics are helpful. Letting $y=Y/N$
and $\Delta y$ denote the increment in $y$ over a time step of size
$1/N$, we find $\mathbb{E}\Delta y = O(1/N)$ while $\mathbb{E}(\Delta
y)^2 = O(1/N^2)$, which means $\var(\Delta y)=O(1/N^2)$. This suggests
that as $N \rightarrow\infty$ we should expect the sample paths of $y$
to approach solutions to the differential equation
%
\begin{equation}
\label{eqe} y' = -r_+y^2 + r_-(1-y).
\end{equation}
Notice the right-hand side is positive at $y=0$, negative at $y=1$ and
strictly decreases with $y$, so there is a unique and globally stable
equilibrium for $y \in[0,1]$, that lies in $(0,1)$. Setting\vspace*{1pt} $y'=0$ and
letting $\alpha= r_+/r_-$ gives the equation $\alpha y^2+y-1=0$ which
has the unique solution $y^*=1/(2\alpha)[-1 + \sqrt{1 + 4\alpha}]$ in
$[0,1]$. Notice that $y^* \sim1 - \alpha$ as $\alpha\rightarrow0^+$
and $y^* \sim1/\sqrt{\alpha}$ as $\alpha\rightarrow\infty$.

\section{Survival analysis}\label{secsurv}
In this section, we analyze $R_0$, $\lambda_c$ and $\Delta(i)$ which
are defined in Section~\ref{secmain}. We begin with $R_0$ defined in
\eqref{eqri}. Define the recruitment probability $p_r =
r_+y^*/(1+r_+y^*) = \mathbb{P}(A\rightarrow E\cup F \cup G)$ which is
the probability of finding a partner before recovering and depends only
on $r_+,r_-$. Define $a = 1+\lambda+r_-$, $b = 2+r_-$ which are the
rates at which the Markov chain of Figure~\ref{figri} jumps away from
states $B$ and $C$, respectively. Also, let
\begin{eqnarray*}
&& \sigma= \sum_{k=0}^{\infty} \biggl(
\frac{\lambda}{a}\frac{2}{b} \biggr)^k = \frac{ab}{ab-2\lambda}.
\end{eqnarray*}
It is easy to check that $ab>2\lambda$. Notice that any path from $A$
to $E\cup F\cup G$ must go to $B$ and then goes around the $B, C$ loop
some number of times before being absorbed at $E,F$ or $G$, and $\sigma
$ accounts for this looping. Summing probabilities over all possible
paths we find
\begin{eqnarray*}
&& \mathbb{P}(A \rightarrow F) = p_r\sigma\frac{r_-}{a}\quad
\mbox{and} \quad\mathbb{P}(A\rightarrow G) = p_r\sigma
\frac{\lambda}{a}\frac{r_-}{b}
\end{eqnarray*}
so we obtain the explicit expression
\begin{eqnarray*}
&& R_0 = p_r\sigma r_-(1 + 2\lambda/b)/a
\end{eqnarray*}
which after re-substituting and a bit of algebra gives
%
\begin{equation}
\label{r0eq0} R_0 = p_rr_-\frac{b+2\lambda}{ab-2\lambda} =
p_rr_-\frac{2+r_-+2\lambda
}{2+3r_-+\lambda r_-+r_-^2}.
\end{equation}

\begin{lemma}\label{r0up}
Fixing $r_+$ and $r_-$, $R_0$ is continuous and increasing with respect
to $\lambda$.
\end{lemma}

\begin{pf}
Continuity is obvious from the formula above. We write $R_0(\lambda)$
and compute the derivative $R_0'(\lambda)$, noting that $p_r$ is fixed.
Letting $c_1=2+r_-$, $c_2=2$, $c_3=2+3r_-+r_-^2$ and $c_4=r_-$,
$R_0(\lambda)=p_rr_-(c_1+c_2\lambda)/(c_3+c_4\lambda)$ so $R_0'(\lambda
) = p_rr_-(c_2c_3-c_1c_4)/(c_3+c_4\lambda)^2$ and $c_2c_3-c_1c_4 =
4+4r_-+r_-^2 > 0$ so $R_0'(\lambda)>0$.
\end{pf}
From this, it follows that for fixed $r_+,r_-$, if $R_0(\lambda)=1$ has
a solution then it is unique and is equal to $\lambda_c$. So, setting
$R_0=1$ gives
%
\begin{equation}
\label{r0eq1} p_rr_-(2+r_-+2\lambda_c) = 2+3r_-+
\lambda_c r_-+r_-^2.
\end{equation}
To get a handle on this equation, we first examine the limit of large
$r_+$, that is, quick formation of partnerships. As noted in Section~\ref{secedge}, $y^* \sim1/\sqrt{\alpha} = \sqrt{r_-}/\sqrt{r_+}$ as $\alpha
=r_+/r_-\rightarrow\infty$, so for fixed $r_-$, $r_+y^* \sim\sqrt
{r_-r_+}\rightarrow\infty$, and so $p_r\rightarrow1$, as
$r_+\rightarrow\infty$. Setting $p_r=1$ in the equation above, after
cancelling like terms and dividing both sides by $r_-$ gives
\begin{eqnarray*}
&& \lambda_c = 1 + 2/r_-
\end{eqnarray*}
for fixed $r_-$, when $r_+=\infty$. For the contact process on a large
complete graph $\lambda_c=1$, so here the only difference is the term
$2/r_-$, which makes it harder for the infection to spread when
partnerships last a long time.

Accounting for $p_r$, we still get a fairly nice expression. From \eqref
{r0eq1}, putting all terms involving $\lambda_c$ on the left and all
other terms on the right gives
\begin{eqnarray*}
&& \lambda_c r_-(2p_r-1) = 2 + (3-2p_r)r_-
+ r_-^2(1-p_r).
\end{eqnarray*}
Letting $\beta= 2p_r-1$ then substituting for $\beta$ and dividing by
$r_-$ gives
%
\begin{equation}
\label{r0eq2} \lambda_c\beta= 2/r_- + (2-\beta) + (1/2)r_-(1-
\beta).
\end{equation}
This equation suggests that we view $\lambda\beta$ as a sort of force
of infection, which makes sense as $\lambda$ is the transmission rate
and $\beta=2p_r-1$ measures the chance of finding a partner before
recovering. Although $\beta$ depends on $r_-$, $-1\leq\beta\leq1$
regardless, so we see from \eqref{r0eq2} that for fixed $\lambda$, if
$r_-$ is either too small or too large, the infection cannot spread.
The reason for this can be understood as follows: if $r_-$ is too
small, partners tend both to recover before breaking up and
transmitting the infection to anyone else, whereas if $r_-$ is too
large, partnerships do not last long enough for transmission to occur.

Using \eqref{r0eq2}, we can now prove the first assertion of Theorem~\ref{thm2}.

\begin{proposition}\label{thm2.1}
For fixed $r_+,r_-$ and $\lambda_c$ given by \eqref{eqlc}, $\lambda
_c<\infty$ if and only if $r_+y^*>1$, if and only if $r_+>1+1/r_-$ and
in this case
\begin{eqnarray*}
&& \lambda_c = \frac{2}{r_-}\frac{2}{(r_+y^*-1)} +
\frac{2}{r_-} + \frac
{4}{r_+y^*-1} + 1 + \frac{r_-}{r_+y^*-1}.
\end{eqnarray*}
\end{proposition}

\begin{pf}
It is easy to check, using the formula $y^* =
(r_-/(2r_+))(-1+(1+4r_+/r_-)^{1/2})$, that $r_+y^*>1$ if and only if
$r_+>1+1/r_-$. Since $\beta\in[-1,1]$, the right-hand side of \eqref
{r0eq2} is positive, so to have a solution it is necessary that $\beta
>0$; dividing by $\beta$ on both sides shows that it is also
sufficient. Then observe that $\beta>0$ if and only if $r_+y^*>1$. To
get the formula for $\lambda_c$, divide by $\beta$ in \eqref{r0eq2} and
observe that $\beta^{-1}-1=2/(r_+y^*-1)$.
\end{pf}

\begin{figure}

\includegraphics{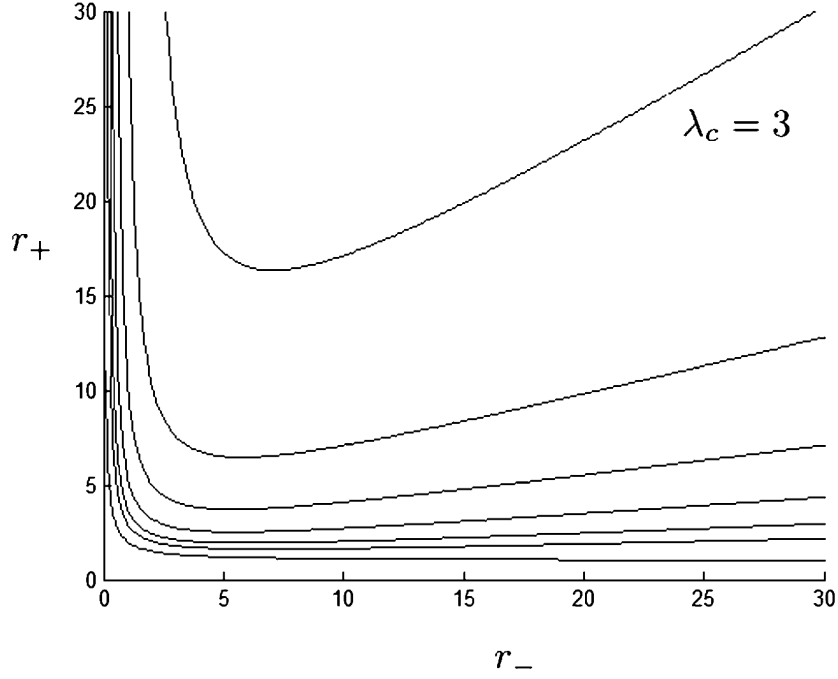}

\caption{Level curves of $\lambda_c$ depicted in the $r_+,r_-$ plane.
Starting from the top curve and going down, $\lambda
_c=3,5,8,13,21,34,\infty$.}
\label{figlc}
\end{figure}

Figure~\ref{figlc} shows level curves of $\lambda_c$ in the $r_+,r_-$
plane. Using the formula for $\lambda_c$, we can see how it scales in
various limits of $r_+,r_-$ and $\alpha$. First, we see what happens
when we speed up and slow down the partnership dynamics. Let $\alpha$
be fixed (and by extension, $y^*$) and let $r_-^*$ denote the unique
value of $r_-$ such that $r_+y^*=1$. We find that:
\begin{itemize}
\item $\lambda_c \downarrow1 + 1/(\alpha y^*)$ as $r_+ \uparrow\infty$ (fast partner
dynamics),
\item $\lambda_c (r_+y^*-1)\downarrow4/r_-^* + 4 + r_-^*$ as $r_+y^*
\downarrow1$ (slow partner dynamics).
\end{itemize}
In particular, in the limit of fast partner dynamics $\lambda_c$
approaches its value for the contact process on a complete graph, plus
a correction for the proportion of available singles. In the slow
limit, that is, as the recruitment probability approaches $1/2$,
$\lambda_c$ diverges like $1/(r_+y^*-1)$, with a proportionality that
itself diverges as $r_-^*$ approaches either $0$ or $\infty$. Now we
fix $r_+>1$ and vary $r_-$. Note that $y^* \downarrow0$ as $r_-
\downarrow0$:
\begin{itemize}
\item as $r_- \uparrow\infty$, $y^* \uparrow1$, $\alpha\downarrow0$
and $\lambda_c/r_- \downarrow1/(r_+-1)$, and
\item as $r_+y^* \downarrow1$, $\lambda_c (r_+y^*-1) \downarrow
4/r_-^* + 4 + r_-^*$.
\end{itemize}
Here, in both limits $\lambda_c$ diverges, in the first case like $r_-$
and in the second case like $1/(r_+y^*-1)$. Finally, we fix $r_-$ and
vary $r_+$, and we find that:
\begin{itemize}
\item as $r_+ \uparrow\infty$, $y^* \sim1/\sqrt{\alpha} = \sqrt
{r^-/r^+}$ and $\lambda_c \rightarrow1 + 2/r_-$, and
\item as $r_+y^* \downarrow 1$, $\lambda_c (r_+y^*-1) \downarrow4/r_-
+ 4 + r_-$.
\end{itemize}
The first limit agrees with the previous large $r_+$ approximation, and
the second limit shows that when $r_+y^*$ is close to 1, $\lambda
(r_+y^*-1)/2 \approx\lambda(r_+y^*-1)/(r_+y^*+1) = \lambda\beta$
behaves like the force of infection and we require again that $r_-$ be
neither too small nor too large in order for the infection to be able
to spread.

We now examine $\Delta(i)$, defined in \eqref{eqdi}.

\begin{lemma}\label{deltalemma}
$\Delta(0)=R_0-1$, and:
\begin{itemize}
\item if $R_0<1$ the equation $\Delta(i)=0$ has no solution $i \in[0,y^*]$,
\item if $R_0=1$ the equation $\Delta(i)=0$ has the unique solution
$i=0$ and
\item if $R_0>1$ the equation $\Delta(i^*)=0$ has a unique solution
$i^* \in(0,y^*)$.
\end{itemize}
\end{lemma}

\begin{pf}
Letting $z = 1+r_+(y^*-i/2)$ we recall the definition:
%
\begin{equation}
\label{deq1} \Delta(i) = p_S\Delta_S +
p_{\mathit{II}}\Delta_{\mathit{II}} + p_{\mathit{SI}}\Delta_{\mathit{SI}}
\end{equation}
with $p_S = 1/z$, $p_{\mathit{SI}} = r_+(y^*-i)/z$, $p_{\mathit{II}} = r_+i/(2z)$, $\Delta
_S=-1$, $\Delta_{\mathit{II}} = -2+\mathbb{P}(C\rightarrow F)+2\mathbb
{P}(C\rightarrow G)$ and $\Delta_{\mathit{SI}} = -1+\mathbb{P}(B\rightarrow
F)+2\mathbb{P}(B\rightarrow G)$, where probabilities are with respect
to the Markov chain in Figure~\ref{figri}.

First, we show $\Delta(0)=R_0-1$. If $i=0$ then $p_S = 1/(1+r_+y^*) =
\mathbb{P}(A\rightarrow D)$, $p_{\mathit{II}}=0$ and $p_{\mathit{SI}} = r_+y^*/(1+r_+y^*)
= \mathbb{P}(A\rightarrow B)$ so
\begin{eqnarray*}
\Delta(0) &=& -\mathbb{P}(A\rightarrow D) + \mathbb{P}(A\rightarrow B) \bigl(-1
+ \mathbb{P}(B\rightarrow F) + 2\mathbb{P}(B\rightarrow G)\bigr)
\\
&=& -\mathbb{P}(A\rightarrow D\cup B) + \mathbb{P}(A\rightarrow F)+2\mathbb{P}(A
\rightarrow G)
\\
&=& -1 + R_0.
\end{eqnarray*}
It is easy to check that $\Delta_{\mathit{II}}\leq0$, so if $\Delta_{\mathit{SI}}\leq0$
then $\Delta(i)<0$ for $i \in[0,y^*]$, since $p_S>0$ and $\Delta_S<0$,
and the other terms are $\leq0$. Since $\partial_i z = -r_+/2$, we find
\begin{eqnarray*}
&& \partial_i p_S = r_+/2z^2 >0\quad
\mbox{and}\quad\partial_i p_{\mathit{II}} = r_+/(2z) +
r_+^2i/\bigl(4z^2\bigr) >0
\end{eqnarray*}
and since $p_{\mathit{SI}}=1-(p_S+p_{\mathit{II}})$, $\partial_i p_{\mathit{SI}}=-\partial_ip_S
-\partial_i p_{\mathit{II}}<0$. If $\Delta_{\mathit{SI}}>0$ it follows that $\partial_i
\Delta(i)<0$ so if $R_0<1$ but $\Delta_{\mathit{SI}}>0$ then $\Delta(i) \leq
\Delta(0)<0$ for $i \in[0,y^*]$. If $R_0 \geq1$, then since $0 \leq
\Delta(0) = p_S\Delta_S + p_{\mathit{SI}}\Delta_{\mathit{SI}}$ and $\Delta_S<0$, it
follows that $\Delta_{\mathit{SI}}>0$ and so $\partial_i \Delta(i)<0$. If
$R_0=1$, then since $\Delta(0)=0$ it follows that $i=0$ is the only
solution in $[0,y^*]$ to the equation $\Delta(i)=0$. If $i=y^*$ then
$p_{\mathit{SI}}=0$ so $\Delta(y^*)<0$, and clearly $\Delta(i)$ is continuous on
$[0,y^*]$. Therefore, if $R_0>1$ then since $\Delta(0)>0$, by the
intermediate value theorem the equation $\Delta(i^*)$ has a solution
$i^* \in(0,y^*)$, and since $\partial_i \Delta(i)<0$ the solution is unique.
\end{pf}

Write $\Delta(i)$ as $\Delta(\lambda,i)$ to emphasize the $\lambda$
dependence. By Lemma~\ref{deltalemma} and since $R_0=1\Leftrightarrow
\lambda=\lambda_c$ and $R_0>1\Leftrightarrow\lambda>\lambda_c$, for
fixed $r_+,r_-$ such that $r_+y^*>1$, we have a function $i^*(\lambda)$
defined for $\lambda\geq\lambda_c$ satisfying $\Delta(\lambda
,i^*(\lambda))=0$ such that $i^*(\lambda_c)=0$ and $i^*(\lambda)>0$ for
$\lambda>\lambda_c$. Next, we see how $i^*$ behaves for $\lambda>\lambda
_c$ near~$\lambda_c$. As usual, $C^1$ means continuously differentiable.

\begin{proposition}\label{thm2.2}
For fixed $r_+,r_-$ such that $r_+y^*>1$, $i^* \sim C(\lambda-\lambda
_c)$ as $\lambda\downarrow\lambda_c$ for some constant $C>0$.
\end{proposition}

\begin{pf}
Clearly, $p_S,p_{\mathit{SI}}$ and $p_{\mathit{II}}$ depend only on $i$ and are $C^1$ in
a neighbourhood of $0$. Also, $\Delta_S$ is fixed and $\Delta_{\mathit{SI}}$ and
$\Delta_{\mathit{II}}$ depend only on $\lambda$ and are rational functions of
$\lambda$ whose range lies in a bounded interval, thus are $C^1$ in a
neighbourhood of $\lambda_c$. Glancing at \eqref{deq1}, this means that
$\Delta(\lambda,i)$ is $C^1$ in a neighbourhood of $(\lambda_c,0)$. If
$\lambda\geq\lambda_c$ then $R_0\geq1$, so as shown in the proof of
Lemma~\ref{deltalemma}, $\partial_i \Delta(\lambda,i)<0$ and in
particular, $\partial_i \Delta(\lambda_c,0)\neq0$. Applying the
implicit function theorem, there is a unique $C^1$ function $i^*(\lambda
)$ defined in a neighbourhood of $\lambda_c$ (and thus coinciding with
the previous definition of $i^*(\lambda)$ when $\lambda\geq\lambda_c$)
satisfying $\Delta(\lambda,i^*(\lambda))=0$, and noting that
$i^*(\lambda_c)=0$,
\begin{eqnarray*}
&& i^*(\lambda) \sim-(\lambda-\lambda_c)\frac{\partial_{\lambda}\Delta
(\lambda_c,0)}{\partial_i \Delta(\lambda_c,0)}
\end{eqnarray*}
as $\lambda\downarrow\lambda_c$. A straightforward Markov chain
coupling argument shows that $\partial_{\lambda}\Delta_{\mathit{SI}},\partial
_{\lambda}\Delta_{\mathit{II}}>0$, which implies $\partial_{\lambda}\Delta
(\lambda,i)>0$. Since $\partial_i\Delta(\lambda,i)<0$, the result follows.
\end{pf}

\section{Mean-field equations}\label{secmf}
A set of differential equations defined below are indispensable to our
analysis of the partner model as they enable a (better and better as
$N$ increases) approximate description of the model, when $N$ is large.
First, we write down the transitions for the variables introduced in
Section~\ref{secmain} that track the total number of singles and pairs
of various types; there are ten such transitions. The existence of
well-defined transitions shows that $(S_t,I_t,\mathit{SS}_t,\mathit{SI}_t,\mathit{II}_t)$ is a
continuous time Markov chain:
\begin{itemize}
\item$I\rightarrow I-1$ and $S\rightarrow S+1$ at rate $I$,
\item$S \rightarrow S-2$ and $\mathit{SS} \rightarrow \mathit{SS}+1$ at rate $(r_+/N)S(S-1)/2$,
\item$S\rightarrow S-1$, $I \rightarrow I-1$ and $\mathit{SI} \rightarrow \mathit{SI}+1$
at rate $(r_+/N)\cdot S \cdot I$,
\item$I \rightarrow I-2$ and $\mathit{II} \rightarrow \mathit{II}+1$ at rate $(r_+/N)I(I-1)/2$,
\item$\mathit{SI} \rightarrow \mathit{SI}-1$ and $\mathit{SS} \rightarrow \mathit{SS}+1$ at rate $\mathit{SI}$,
\item$\mathit{II} \rightarrow \mathit{II}-1$ and $\mathit{SI} \rightarrow \mathit{SI}+1$ at rate $2\mathit{II}$,
\item$\mathit{SI} \rightarrow \mathit{SI}-1$ and $\mathit{II} \rightarrow \mathit{II}+1$ at rate $\lambda \mathit{SI}$,
\item$\mathit{SS} \rightarrow \mathit{SS}-1$ and $S \rightarrow S+2$ at rate $r_-\mathit{SS}$,
\item$\mathit{SI} \rightarrow \mathit{SI}-1$, $S \rightarrow S+1$ and $I\rightarrow I+1$
at rate $r_-\mathit{SI}$, and
\item$\mathit{II} \rightarrow \mathit{II}-1$ and $I \rightarrow I+2$ at rate $r_-\mathit{II}$.
\end{itemize}
Focusing now on the rescaled quantities
$(s_t,i_t,ss_t,si_t,ii_t)=(S_t,I_t,\mathit{SS}_t,\mathit{SI}_t,\break 
\mathit{II}_t)/N$ and noting the
relation $s_t+i_t+2(ss_t+si_t+ii_t)=1$, we shall ignore $ss_t$ since it
plays no role in the calculations that follow. Also, it will be
convenient to use $y_t:=s_t+i_t$ instead of $s_t$. Doing so, the above
transitions become:
\begin{itemize}
\item$i\rightarrow i-1/N$ at rate $iN$,
\item$y \rightarrow y-2/N$ at rate $[r_+(y-i)(y-i-1/N)/2]N =
[r_+(y-i)^2/2]N - r_+(y-i)/2$,
\item$y \rightarrow y-2/N$, $i \rightarrow i-1/N$ and $si \rightarrow
si+1/N$ at rate $r_+(y-i)iN$,
\item$y \rightarrow y-2/N$, $i \rightarrow i-2/N$ and $ii \rightarrow
ii+1/N$ at rate $[r_+i(i-1/N)/2]N = (r_+i^2/2)N - r_+i/2$,
\item$si \rightarrow si-1/N$ at rate $siN$,
\item$ii \rightarrow ii-1/N$ and $si \rightarrow si+1/N$ at rate $2iiN$,
\item$si \rightarrow si-1/N$ and $ii \rightarrow ii+1/N$ at rate
$\lambda siN$,
\item$y\rightarrow y+2/N$ at rate $[r_-((1-y)/2-(si+ii))]N$,
\item$si \rightarrow si-1/N$, $y \rightarrow y+2/N$ and $i\rightarrow
i+1/N$ at rate $r_-siN$, and
\item$ii \rightarrow ii-1/N$, $y\rightarrow y+2/N$ and $i \rightarrow
i+2/N$ at rate $r_-iiN$.
\end{itemize}
As we did for $y_t$ in Section~\ref{secedge}, we derive some
differential equations that approximate the evolution of
$(y_t,i_t,si_t,ii_t)$; since we already have an equation for $y_t$ we
focus on $i_t,si_t,ii_t$. We have
\begin{eqnarray*}
\frac{d}{dt}\mathbb{E}(i_t \vert i_t=i)
&=& -\bigl(1+r_+(y-i)+2r_+(i-1/N)/2\bigr)i + r_-(si+2ii),
\\
\frac{d}{dt}\mathbb{E}(si_t \vert si_t=si) &=&
r_+(y-i)i+2ii-(1+\lambda +r_-)si,
\\
\frac{d}{dt}\mathbb{E}(ii_t \vert ii_t=ii) &=&
r_+i(i-1/N)/2 +\lambda si - (2+r_-)ii
\end{eqnarray*}
and as before, in a time step of size $1/N$ the increment in each
variable has expected value $O(1/N)$ while its square has expected
value $O(1/N^2)$. Adding in the $y'$ equation \eqref{eqe}, this
suggests again that in the limit as $N\rightarrow\infty$ we should
expect the sample paths of $(y_t,i_t,si_t,ii_t)$ to approach solutions
to the \emph{mean-field equations}
%
\begin{eqnarray}
y' &=& -r_+y^2 + r_-(1-y),
\nonumber
\\
i' &=& -(1+r_+y)i + r_-(si+2ii),
\nonumber
\\[-8pt]
\label{mfeq}
\\[-8pt]
\nonumber
si' &=& r_+(y-i)i -(1+\lambda+r_-)si + 2ii,
\\
\nonumber
ii' &=& r_+i^2/2 + \lambda si -(2+r_-)ii.
\end{eqnarray}
It is sometimes convenient to replace $si$ with $ip := si+ii$, where
the $ip$ stands for ``infected partnership''. Since $si=ip-ii$, both
forms lead to the same solutions. After the change of variables, we have
%
\begin{eqnarray}
y' &=& -r_+y^2 + r_-(1-y),
\nonumber
\\
i' &=& -(1+r_+y)i + r_-(ip+ii),
\nonumber
\\[-8pt]
\label{mfeq2}
\\[-8pt]
\nonumber
ip' &=& r_+(y-i/2)i-(1+r_-)ip+ii,
\\
\nonumber
ii' &=& r_+i^2/2 + \lambda ip -(2+r_-+
\lambda)ii.
\end{eqnarray}
We will often use the shorthand $u'=F(u)$ for the MFE \eqref{mfeq} or
\eqref{mfeq2}, where $u\in\mathbb{R}^4$. In both cases the MFE have
the form $y'=f(y),u'=G(y,u)$, where \mbox{$u\in\mathbb{R}^3$}, that is, the
$y$ dynamics does not depend on the other 3 variables, but it does
influence them; systems of this form are often referred to as \emph
{skew product}. The next three results have natural analogues for the
stochastic model, and in fact the analogue of Lemma~\ref{mfmt} shows up
in Section~\ref{secmacro} as Lemma~\ref{pmmt}. First, we show the
domain of interest is an invariant set.

\begin{lemma}\label{mfatt}
The following set is invariant for the MFE:
\begin{eqnarray*}
&&\Lambda:= \bigl\{(y,i,ip,ii) \in\mathbb{R}^4_+:i\leq y \leq1,ii
\leq ip \leq(1-y)/2\bigr\}.
\end{eqnarray*}
\end{lemma}

\begin{pf}
We examine the boundary and use the form \eqref{mfeq2} of the MFE. If
$y=0$ then $y'>0$ and if $y=1$ then $y'<0$, so $[0,1]$ is invariant for
$y$. Let $u=(i,ip,ii)$. If $u=(0,0,0)$, then $u'=(0,0,0)$, so $(0,0,0)$
is invariant for $u$. If $u\neq(0,0,0)$ and $u_j=0$ for coordinate $j$,\vspace*{-1pt}
then $u_j'>0$ (note for $ip'$ that since $i\leq y$, if $i>0$ then
$y-i/2>0$). If $i=y \neq0$, then since $ip+ii\leq(1-y)$, $i' \leq-y
- r_+y^2 + r_-(1-y) =-y+y'< y'$. If $i=y=0$ then $i' \leq-y+y'=y'$ and
since $y'>0$, $i'' \leq-y' + y'' < y''$. For the remainder, we may
assume $i<y$. If $ii=ip \neq0$, then $ii' = r_+i^2/2 -(2+r_-)ip \leq
r_+(y-i/2)i -(2+r_-)ip < ip'$ while if $ii=ip=0$ then we may assume
$i>0$ in which case $ii' = r_+i^2/2 < r_+(y-i/2)i = ip'$.
\end{pf}
Written in the form \eqref{mfeq2}, the MFE have a useful monotonicity
property which is described in the following lemma.

\begin{lemma}\label{mfmt}
Let $(y(t),u(t))$ and $(y(t),v(t))$ be solutions to the MFE written in
$(y,i,ip,ii)$ coordinates, and say that $u\leq v \Leftrightarrow
u_j\leq v_j \  \forall j \in\{1,2,3\}$. If $u(0)\leq v(0)$, then
$u(t)\leq v(t)$ for $t>0$.
\end{lemma}

\begin{pf}
Since trajectories are continuous it suffices to check that if $u\leq
v$, $u\neq v$ and $u_j=v_j$ then $u_j' < v_j'$. Referring to \eqref
{mfeq2}, $i'$ increases with $ip$ and $ii$, $ip'$ increases with $i$
and $ii$ [note $\partial_i(y-i/2)i = y-i$ and $i\leq y$] and $ii'$
increases with $i$ and $ip$.
\end{pf}
For what follows, we set $y=y^*$ in which case the MFE are
three-dimensional. Since $\Lambda$ is invariant,
\begin{eqnarray*}
&& \Lambda^* := \bigl\{(y,u) \in\Lambda:y=y^*\bigr\}
\end{eqnarray*}
is also invariant. Since $\Lambda^* \cong\{(i,ip,ii)\in\mathbb
{R}^3_+:i\leq y^*, ii\leq ip \leq(1-y^*)/2\}$ is three-dimensional,
elements of $\Lambda^*$ are usually written as a three-vector in either
$(i,si,ii)$ or $(i,ip,ii)$ coordinates.

\begin{lemma}\label{mfthmsuff}
Say that $u=(i,ip,ii)$ is \emph{increasing} if $u_j'>0$ in each
coordinate. For the MFE with $y=y^*$ and any solution $u(t)$:
\begin{itemize}
\item if $(0,0,0)$ is the only equilibrium then $u(t)\rightarrow
(0,0,0)$ as $t\rightarrow\infty$, and
\item if there is a unique equilibrium $u^*\neq(0,0,0)$ and a sequence
of nonzero increasing states tending to $(0,0,0)$, then for $u(0)\neq
(0,0,0)$, $u(t)\rightarrow u^*$ as $t\rightarrow\infty$.
\end{itemize}
\end{lemma}

\begin{pf}
Defining $\overline{u}:=(y^*,(1-y^*)/2,(1-y^*)/2)$, $\overline{u}\geq
v$ for all $v\in\Lambda^*$, so letting $\overline{u}(t)$ be the
solution to the MFE with $\overline{u}(0)=\overline{u}$, for $s\geq0$,
$\overline{u}(0)\geq\overline{u}(s)$. Since $y=y^*$, by monotonicity
(Lemma~\ref{mfmt}) $\overline{u}(t)\geq\overline{u}(t+s)$ for $t>0$,
so $\overline{u}(t)$ is nonincreasing in $t$. Since $\Lambda^*$ is
compact, $\lim_{t\rightarrow\infty}\overline{u}(t)$ exists and by
continuity of the MFE is an equilibrium. If $(0,0,0)$ is the only
equilibrium, then since $\overline{u}(t)\geq(0,0,0)$, $\overline
{u}(t)\rightarrow(0,0,0)$ as $t\rightarrow\infty$, so for any solution
$v(t)$, since $\overline{u}(0)\geq v(0)$, $\overline{u}(t)\geq v(t)$
for $t>0$, and since $v(t)\geq(0,0,0)$, $v(t)\rightarrow(0,0,0)$.

If $u(0)$ is increasing, then $u(0)\neq(0,0,0)$ and by continuity of
the MFE there is $\varepsilon>0$ so that $u(s) \geq u(0)$ for $0 \leq s
\leq\varepsilon$. By monotonicity $u(t+s) \geq u(t)$ for $0 \leq s \leq
\varepsilon$ and if $(k-1)\varepsilon\leq s \leq k\varepsilon$, by iterating at
most $k$ times $u(t+s)\geq u(t)$, so $u(t)$ is increasing for all time.
As in the previous case, $\lim_{t\rightarrow\infty}u(t)$ exists and is
an equilibrium which in this case is not $(0,0,0)$. If there is a
unique equilibrium $u^* \neq(0,0,0)$, and if for any nonzero solution
$v(t)$ there is $T>0$ so that $v(T)\geq u$ for some increasing $u$,
then setting $u(T)=u$, since $\overline{u}(t)\geq v(t) \geq u(t)$ for
$t \geq T$ and $\lim_{t\rightarrow\infty}\overline{u}(t)=\lim_{t\rightarrow\infty}u(t) = u^*$ it follows that $\lim_{t\rightarrow
\infty}v_t = u^*$. If $v(0)\neq(0,0,0)$, then for $t>0$, $v_j(t)>0$ in
each coordinate $j$; this follows from the fact that for $j=1,2,3$,
$v_j' \geq-Cv_j$ for some $C$, and if $v_j=0$ but $v_k>0$ for some
$k\neq j$ then $v_j'>0$. Thus, fixing $T>0$, if $v(0)\neq(0,0,0)$ then
since $\varepsilon:=\min_j v_j(T)>0$, if there is a sequence of increasing
states tending to $(0,0,0)$ there is an increasing state $u$ with $\max_j u_j \leq\varepsilon$, and thus $v(T)\geq u$, as desired.
\end{pf}
As the next result shows, on $\Lambda^*$ the MFE have a simple dynamics
with a bifurcation at $R_0=1$. Since we refer back to quantities from
Section~\ref{secsurv}, in this proof we mostly use $(i,si,ii)$ coordinates.

\begin{theorem}\label{mfthm}
For the MFE:
\begin{itemize}
\item if $R_0 \leq1$ there is the unique equilibrium $(0,0,0)$ which
is attracting on $\Lambda^*$ and
\item if $R_0>1$ there is a unique positive equilibrium
$(i^*,s^*,ii^*)$ satisfying $\Delta(i^*)=0$ which is attracting on
$\Lambda^*\setminus\{(0,0,0)\}$.
\end{itemize}
\end{theorem}

\begin{pf}
By Lemma~\ref{mfthmsuff} it is enough to show that if $R_0\leq1$ then
$(0,0,0)$ is the only equilibrium, and that if $R_0>1$ there is a
unique equilibrium $(i^*,si^*,ii^*)\neq(0,0,0)$ satisfying $\Delta
(i^*)=0$, and a sequence of increasing states converging to $(0,0,0)$.
Treating $si,ii$ as a separate system with input function $i$, we have
the nonhomogenous linear system
\begin{eqnarray*}
&& \pmatrix{si'
\cr
ii'} = %
\pmatrix{-a & 2
\cr
\lambda& -b} \pmatrix{si
\cr
ii } %
+ r_+i
\pmatrix{\bigl(y^*-i\bigr)
\cr
i/2 } %
\end{eqnarray*}
or, in matrix form, $v' = Kv + Li$, with $v = (si,ii)^{\top}$, $K =
{-a \quad 2\,\, \choose  \,\,\lambda\!\!\quad -b}$ and $L = r_+((y^*-i),i/2)^{\top}$, whose solution is given by
%
\begin{equation}
\label{subsyssol} v(t) = \Phi(t)v(0) + \int_0^t
\Phi(t-s)L(s)i(s)\,ds,
\end{equation}
where $\Phi(t) = \exp(Kt)$ is the solution of the associated homogenous
system---note that $\Phi(t)$ is the restriction of the transition
semigroup for the continuous-time Markov chain from Figure~\ref{figri}
to the states $B$ and $C$. Substituting the solution for the $si,ii$
system into the equation for $i$, we have
%
\begin{equation}
\label{ieq}\hspace*{6pt} i'(t) = -\bigl(1+r_+y^*\bigr)i(t) + r_-(1,2) \biggl[
\Phi(t)v(0) + \int_0^t \Phi (t-s)L(s)i(s)\,ds
\biggr],\hspace*{-24pt}
\end{equation}
where $(1,2)$ is a row vector that multiplies the column vector in the
square brackets. This equation depends only on $i$, the initial values
$v(0) = (si(0),ii(0))^{\top}$ and the solution matrix $\Phi(t)$.

Linearizing \eqref{ieq} around $(i,si,ii)=(0,0,0)$ and using the ansatz
$i(t)=\exp(\mu t)$, we obtain
\begin{eqnarray*}
&&\mu e^{\mu t} = -\bigl(1+r_+y^*\bigr)e^{\mu t} + r_-(1,2)
\biggl[\Phi(t)v_0 + \int_0^t
\Phi(t-s)e^{\mu s}\,ds \biggr]L_0,
\end{eqnarray*}
where $L_0 = r_+(y^*,0)^{\top}$, and using $\Phi(t)=\exp(Kt)$ the
integral in the square brackets is
\begin{eqnarray*}
&&\!\! e^{Kt}\int_0^t e^{(\mu I-K)s}
\,ds = e^{Kt}(\mu I-K)^{-1}\bigl(e^{(\mu I-K)t} - I\bigr) =
(\mu I-K)^{-1}\bigl(e^{\mu t} - e^{Kt}\bigr),
\end{eqnarray*}
where $I$ is the identity matrix. Letting $t \rightarrow\infty$ and
noting $\Phi(t) = e^{Kt} \rightarrow0$ since $K$ is a stable matrix,
we obtain the eigenvalue equation
\begin{eqnarray*}
&&\mu= -\bigl(1+r_+y^*\bigr) + r_-(1,2) (\mu I-K)^{-1}L_0
\end{eqnarray*}
which, expanding, is
\begin{eqnarray*}
\label{mueq} && \mu= -\bigl(1+r_+y^*\bigr) + r_-\frac{\mu+ b + 2\lambda}{(\mu+b)(\mu
+a)-2\lambda}r_+y^*
\end{eqnarray*}
and setting $\mu=0$ gives the equation
\begin{eqnarray*}
&& 1 = \frac{r_+y^*}{1+r_+y^*}\frac{r_-}{ab-2\lambda}(b+2\lambda)
\end{eqnarray*}
which, comparing to \eqref{r0eq0}, is exactly $R_0=1$. Recalling that
$ab-2\lambda>0$,
\begin{eqnarray*}
&&\frac{d}{d\mu} \biggl( \frac{\mu+b+2\lambda}{(\mu+b)(\mu+a)-2\lambda
} \biggr)\\
 &&\qquad= \frac{(\mu+b)(\mu+a)-2\lambda- (\mu+b+2\lambda)(2\mu+ b
+ a)}{[(\mu+b)(\mu+a)-2\lambda]^2}
\\
&&\qquad= \frac{-2\lambda-[(\mu+b)^2+2\lambda(2\mu+b+a)]}{[(\mu+b)(\mu
+a)-2\lambda]^2}
\end{eqnarray*}
is negative when $\mu\geq0$. Setting $\mu=0$ in \eqref{mueq}, the
right-hand side is positive if $R_0>1$, so since both sides are
continuous in $\mu$, the left-hand side is equal to $0$ at $\mu=0$ and
increases unboundedly as $\mu$ increases and the right-hand side
decreases with $\mu$ it follows that \eqref{mueq} has a positive
solution $\mu>0$ when $R_0>1$.

To obtain the increasing states mentioned in Lemma~\ref{mfthmsuff}, we
show that for \mbox{$R_0>1$} the unstable eigenvector of the linearized system
near $(0,0,0)$ is strictly positive when viewed in $(i,ip,ii)$
coordinates; we can then take for the initial states small multiples of
the eigenvector. To show the eigenvector is strictly positive,
linearize \eqref{subsyssol} around $(i,si,ii)=(0,0,0)$ with input $i(t)
=\exp(\mu t)$, substitute the solution form $v(t) = v \exp(\mu t)$ and
let $t\rightarrow\infty$ to obtain $v = (\mu I-K)^{-1}L_0$ which has
positive entries, which implies that in $(ip,ii)$ coordinates it also
has positive entries.

It remains to look for nonzero equilibria. Focusing again on \eqref
{ieq}, as our steady state assumption we suppose the system was started
in the distant past and has remained in equilibrium up to the present
time. Since\vspace*{1pt} $\Phi(t)\rightarrow0$ as $t \rightarrow\infty$ we ignore
$\Phi(t)v(0)$, and letting\vspace*{1pt} $\Phi_{\infty} = \int_0^{\infty}\Phi(s)\,ds =
-K^{-1}$, $\int_0^t \Phi(t-s)L(s)i(s)\,ds$ becomes $\smash{\int_{-\infty
}^t \Phi(t-s)L^{\dagger}i^{\dagger}\,ds} = \Phi_{\infty}L^{\dagger
}i^{\dagger}$ where $L^{\dagger} = r_+((y^*-i^{\dagger}),i^{\dagger
}/2)^{\top}$ and $i^{\dagger}$ are the equilibrium values, and we obtain
\begin{eqnarray*}
&& \bigl(1+r_+y^*\bigr) = r_-(1,2)\Phi_{\infty}L^{\dagger}.
\end{eqnarray*}
Notice that $r_-(1,2)\Phi_{\infty}$ returns the expected number of
infectious singles that result from an $\mathit{SI}$ or an $\mathit{II}$ partnership upon
breakup, so we have $r_-(1,2)\Phi_{\infty} = (1+\Delta_{\mathit{SI}},2+\Delta
_{\mathit{II}})$ and
\begin{eqnarray*}
\bigl(1+r_+y^*\bigr) &=& r_+\bigl[\bigl(y^*-i^{\dagger}\bigr) (1+
\Delta_{\mathit{SI}}) + \bigl(i^{\dagger
}/2\bigr) (2+\Delta_{\mathit{II}})
\bigr]
\\
&=& r_+y^* + r_+\bigl[\bigl(y^*-i^{\dagger}\bigr)\Delta_{\mathit{SI}} +
\bigl(i^{\dagger}/2\bigr)\Delta_{\mathit{II}}\bigr]
\end{eqnarray*}
and subtracting $r_+y^*$, $1=r_+(y^*-i^{\dagger})\Delta_{\mathit{SI}} +
r_+(i^{\dagger}/2)\Delta_{\mathit{II}}$ which comparing with~\eqref{deq1} is
exactly the equation $\Delta(i^{\dagger})=0$, as desired. By Lemma~\ref
{deltalemma}, we have the unique solution $i^{\dagger} = i^*$ if
$R_0>1$, and there is no positive solution when $R_0\leq1$. Using the
steady state assumption and \eqref{subsyssol} gives $(si^{\dagger
},ii^{\dagger})=\Phi_{\infty}L^{\dagger}i^{\dagger}$, that is,
$si^{\dagger},ii^{\dagger}$ are uniquely determined by $i^{\dagger}$.
This proves uniqueness of the nonzero equilibrium when $R_0>1$ and
uniqueness of $(0,0,0)$ as an equilibrium when $R_0 \leq1$.
\end{pf}

\begin{remark}\label{r0remark}
Setting $y=y^*$ in \eqref{mfeq} and writing the remaining equations in
matrix form, we have $u' = Au$ with $u = (i,si,ii)^{\top}$ and
\begin{eqnarray*}
A = %
\pmatrix{ -\bigl(1+r_+y^*\bigr) & r_- & 2r_-
\cr
r_+\bigl(y^*-i
\bigr) & -a & 2
\cr
r_+i & \lambda& -b } %
\end{eqnarray*}
that depends on $u$. Using the technique of \cite{watm}, if we evaluate
$A$ at $i=0$ and write it as $F-V$ with
\begin{eqnarray*}
F = %
\pmatrix{ 0 & 0 & 0
\cr
r_+y^* & 0 & 0
\cr
0 & 0 & 0 },\qquad V =
\pmatrix{ \bigl(1+r_+y^*\bigr) & -r_- & -2r_-
\cr
0 & a & -2
\cr
0 & -
\lambda & b }
\end{eqnarray*}
and define $R_0 = \rho(FV^{-1})$ where $\rho$ is the spectral radius,
then it can be verified that this definition of $R_0$ coincides with
the one given in \eqref{eqri}. Then, according to Theorem~2 of \cite
{watm}, $R_0<1$ implies $(0,0,0)$ is locally asymptotically stable,
while $R_0>1$ implies it is unstable.
\end{remark}

\section{Approximation by the mean-field equations}\label{secstoch}
In this section, we show how to approximate the sample paths of
$(y_t,i_t,si_t,ii_t)$ with solutions to the MFE~\eqref{mfeq}, and use
this to get some control on $y_t$. Unless otherwise noted, for a
vector, $|\cdot|$ denotes the $\ell^{\infty}$ norm, that is, $|u| = \max_i |u_i|$. We begin with a useful definition.

\begin{definition}\label{defwhp}
An event $A$ depending on a parameter $n$ is said to hold \emph{with
high probability} or w.h.p. in $n$ if there exists $\gamma>0$ and $n_0$ so
that $\mathbb{P}(A) \geq1-e^{-\gamma n}$ when $n\geq n_0$.
\end{definition}

When possible, probability estimates are given more or less explicitly,
but we will occasionally use this definition to reduce clutter,
especially in Section~\ref{secmacro}. We begin with a well-known large
deviations result for Poisson random variables; since it is not hard to
prove, we supply the proof. For a reference to large deviations theory,
see Section~1.9 in \cite{prob}.

\begin{lemma}\label{chern}
Let $X$ be Poisson distributed with mean $\mu$, then
\begin{eqnarray*}
\mathbb{P}\bigl(X>(1+\delta) \mu\bigr) &\leq& e^{-\delta^2\mu/4}\qquad \mbox{for }
0<\delta\leq1/2,
\\
\mathbb{P}\bigl(X<(1-\delta)\mu\bigr) &\leq& e^{-\delta^2\mu/2}\qquad \mbox{for }
\delta>0.
\end{eqnarray*}
\end{lemma}

\begin{pf}
We deal separately with $X > (1+\delta)\mu$ and $X<(1-\delta)\mu$. For
$t>0$ and using Markov's inequality we have
\begin{eqnarray*}
&& \mathbb{P}\bigl(X>(1+\delta)\mu\bigr) = \mathbb{P}\bigl(e^{tX}>e^{(1+\delta)t\mu}
\bigr) \leq\mathbb{E}e^{tX}e^{-(1+\delta)t\mu}.
\end{eqnarray*}
Notice that
\begin{eqnarray*}
&& \mathbb{E}e^{tX} = \sum_{k\geq0}e^{tk}e^{-\mu}
\frac{\mu^k}{k!} = e^{-\mu}\sum_{k\geq0}
\frac{(e^t\mu)^k}{k!} = e^{-\mu}e^{e^t\mu} = \exp \bigl(
\bigl(e^t-1\bigr)\mu\bigr)
\end{eqnarray*}
so $\mathbb{E}e^{tX}e^{-(1+\delta)t\mu}=\exp(\mu(e^t-1-(1+\delta)t))$.
Minimizing $e^t-1-(1+\delta)t$ gives $t=\log(1+\delta)$, and thus
$(1+\delta)-1-(1+\delta)\log(1+\delta)=\delta-(1+\delta)\log(1+\delta
)$. Since $\log(1+\delta) \geq\delta-\delta^2/2$ this is at most
$\delta-(1+\delta)(\delta-\delta^2/2) = -\delta^2/2+\delta^3/2$ which
is $\leq-\delta^2/4$ for $0<\delta\leq1/2$.

For the other direction we take a similar approach. For $t>0$ and
using Markov's inequality we have
\begin{eqnarray*}
&& \mathbb{P}\bigl(X<(1-\delta)\mu\bigr) = \mathbb{P}\bigl(e^{-tX}>e^{-(1-\delta)t\mu}
\bigr) \leq\mathbb{E}e^{-tX}e^{(1-\delta)t\mu}
\end{eqnarray*}
and using $\mathbb{E}e^{-tX} = \exp((e^{-t}-1)\mu)$ the right-hand side
above is $\exp(\mu(e^{-t}-1+(1-\delta)t))$. Minimizing
$e^{-t}-1+(1-\delta)t$ gives $-t=\log(1-\delta)$, and thus $(1-\delta
)-1-(1-\delta)\log(1-\delta)=-\delta-(1-\delta)\log(1-\delta)$. Since
$\log(1-\delta) \geq-\delta-\delta^2/2$, this is at most $-\delta
+(1-\delta)(\delta+\delta^2/2) = -\delta^2/2-\delta^3/2 \leq-\delta^2/2$.
\end{pf}
For the next three results, we use the notation $u_t =
(y_t,i_,si_t,ii_t)$. First, we give an a priori bound on the change in
$u_t$ over a short period of time.

\begin{lemma}\label{apriori}
Let $u_t=(y_t,i_t,si_t,ii_t)$. There are constants $C,\gamma>0$ so that
for all $h>0$ and fixed $t$,
\begin{eqnarray*}
&& \mathbb{P}\Bigl(\sup_{t\leq s \leq t+h}|u_s-u_t|
\leq Ch\Bigr) \geq1-e^{-\gamma Nh}.
\end{eqnarray*}
\end{lemma}

\begin{pf}
Looking to the transitions listed in Section~\ref{secmf}, jumps in
$u_t$ are of size $\leq2/N$ and occur at total rate $\leq MN$ for some
$M>0$ that depends only on parameters. Thus, in a time step $h>0$ the
number of events affecting $u_t$ is stochastically bounded above by a
Poisson random variable $X$ with mean $MNh$, so if $X \leq x$ then
$|u_s-u_t|\leq2x/N$ for all $s \in[t,t+h]$. By Lemma~\ref{chern},
$\mathbb{P}(X>(1+\delta)MNh) \leq e^{-\delta^2MNh/4}$ for $0<\delta\leq
1/2$. Taking $\delta=1/4$ and $C = 2(1+\delta)M$, $\gamma= \delta
^2M/4$ completes the proof.
\end{pf}
Let $u'=F(u)$ denote the MFE \eqref{mfeq}. As $N$ becomes large, for
small $h>0$ we expect that with probability tending to 1, $u_{t+h} =
u_t + hF(u_t) + o(h)$. Using Lemma~\ref{apriori} and re-using the
estimate from Lemma~\ref{chern} we obtain a quantitative bound on the remainder.

\begin{lemma}\label{dapriori}
Let $u_t=(y_t,i_t,si_t,ii_t)$. For each $\varepsilon>0$ there are
constants $C,\gamma>0$ so that for small enough $h>0$,
\begin{eqnarray*}
&& \mathbb{P}\bigl(\bigl|u_{t+h}-u_t-hF(u)\bigr|\leq\varepsilon h
\bigr) \geq1-Ce^{-\gamma Nh}.
\end{eqnarray*}
\end{lemma}

\begin{pf}
Let $Q_j(u)$, $j=1,\ldots,10$, denote the transition rates of the ten
transitions introduced in Section~\ref{secmf}, as a function of $u$,
and let $X_j(t,h)$ denote the number of type $j$ transitions occurring
in the time interval $[t,t+h]$. For each $j$, $Q_j(u)=Nq_j(u) + R_j(u)$
where $q_j(u)$ is a quadratic function of $u$ and $R_j(u)$ is a
remainder that satisfies $|R_j(u)|\leq M$ for some $M>0$ and all $u \in
[0,1]^4$. It is easily verified that if $u_t=u$ and $X_j(t,h)=Nq_j(u)h$
for each $j$ then $u_{t+h}=u + hF(u)$. Since each transition changes
$u$ by at most $2/N$, it is therefore enough to show that there are
constants $C,\gamma>0$ so that for each $j$, small enough $h>0$, and
all $u$,
\begin{eqnarray*}
&& \mathbb{P}\bigl(\bigl|X_j(t,h)-Nq_j(u)h\bigr|\leq\varepsilon
Nh/20 \vert u_t=u\bigr)\geq 1-Ce^{-\gamma Nh}.
\end{eqnarray*}
Since the domain of $q_j(u)$ is a subset of $[0,1]^4$, and thus bounded
it follows that $q_j$ is bounded and Lipschitz continuous, that is, for
some $L>0$ and all $v,u$ in the domain of $q_j$, $q_j(u)\leq L$ and
$|q_j(v)-q_j(u)|\leq L|v-u|$, and in particular, $|Q_j(v)-Q_j(u)| \leq
NL|v-u| + 2M$; for what follows, take $L\geq\varepsilon$. Let $A(t,h)$ be
the event
\begin{eqnarray*}
&& \Bigl\{\sup_{t\leq s \leq t+h}|u_s-u_t|\leq
C_1h\Bigr\},
\end{eqnarray*}
from Lemma~\ref{apriori}, then on the event $\{u_t=u\}\cap A(t,h)$,
\begin{eqnarray*}
 \sup_{t\leq s \leq t+h}\bigl|Q_j(u_s)-Nq_j(u)\bigr|
&\leq & \sup_{t\leq s \leq
t+h}\bigl|Q_j(u_s)-Q_j(u)\bigr|
+ \bigl|Q_j(u)-Nq_j(u)\bigr|\\
& \leq &  N(LC_1h + 3M/N).
\end{eqnarray*}
For ease of notation, let $q=q_j(u)$ and let $r=LC_1h+3M/N$, and note
that $r \rightarrow0$ as $\max(h,1/N) \rightarrow0$. Then, on $\{
u_t=u\}\cap A(t,h)$, $X_j(t,h)$ is stochastically bounded above and
below respectively by Poisson random variables with means $Nh(q+r)$ and
$Nh(q-r)$, so from Lemma~\ref{chern} it follows that for $0<\delta\leq1/2$,
%
\begin{eqnarray}
&& \mathbb{P}\bigl(\bigl\{\bigl|X_j(t,h) - Nhq\bigr|
\leq Nh
\bigl(q\delta+ r(1+\delta)\bigr)\bigr\} \cap\{ u_t=u\}\cap A(t,h)
\bigr)
\nonumber
\\[-8pt]
\label{derivest}
\\[-8pt]
\nonumber
&&\qquad \geq1-2e^{-Nh(q-r)\delta^2/4}.
\end{eqnarray}
Recalling that $q \leq L$, let $h,\delta,1/N>0$ be chosen small enough
that $L\delta+ r(1+\delta)\leq\varepsilon/20$, then $Nh(q\delta+
r(1+\delta)) \leq\varepsilon Nh/20$. To bound the probability uniformly
in $q$, we split into two cases according as $q\geq q\delta+ r(1+\delta
)$ or not, that is, as $q \geq r(1+\delta)/(1-\delta)$ or not. If
$q\geq r(1+\delta)/(1-\delta)$ then letting $\gamma_1 = r[(1+\delta
)/(1-\delta)-1]\delta^2/4$ which is $>0$ it follows that $Nh(q-r)\delta
^2/4 \geq\gamma_1 Nh$. If $q<q\delta+r(1+\delta)$ the lower bound on
$X_j(t,h) - Nhq$ is trivial and so in that case
\begin{eqnarray*}
&& \mathbb{P}\bigl(\bigl\{\bigl|X_j(t,h)-Nhq\bigr| \leq Nh\bigl(q\delta+ r(1+
\delta)\bigr)\bigr\}\cap\{ u_t=u\}\cap A(t,h)\bigr)
\\
&&\qquad\geq 1-e^{-Nh(q+r)\delta^2/4}.
\end{eqnarray*}
Letting $\gamma_2 = r\delta^2/4$ which is $>0$ it follows that
$Nh(q+r)\delta^2/4 \geq\gamma_2 Nh$. Letting $\gamma_3$ be such that
$\mathbb{P}(A(t,h))\geq1-e^{-\gamma_3Nh}$ and letting $\gamma= \min
(\gamma_1,\gamma_2,\gamma_3)$ and $C=3$ completes the proof.
\end{pf}
Using the above estimate, we obtain finite-time control on the
evolution of $u_t$, as $N$ becomes large.

\begin{proposition}\label{mfest}
Let $u_t=(y_t,i_t,si_t,ii_t)$. For each $\varepsilon,T>0$ there are
constants $\delta,C,\gamma>0$ so that from any initial condition $u_0$
and any solution $u(t)$ to the MFE \eqref{mfeq} satisfying
$|u_0-u(0)|\leq\delta$,
\begin{eqnarray*}
&& \mathbb{P}\Bigl(\sup_{0 \leq t \leq T}\bigl|u_t-u(t)\bigr|\leq
\varepsilon\Bigr) \geq 1-Ce^{-\gamma N}.
\end{eqnarray*}
\end{proposition}

\begin{pf}
The proof is analogous to the proof in numerical analysis that the
Euler method is $O(h)$ accurate. Fix $h=T/M$ for integer $M$ and define
events $A_1,\ldots,A_m$ as follows: $A_1=B_1\cap D_1$ and given $A_{j-1}$,
$A_j=A_{j-1} \cap B_j \cap D_j$ where
\begin{eqnarray*}
&& B_j = \Bigl\{\sup_{h(j-1) \leq t \leq hj}|u_t -
u_{hj}| \leq C_1h\Bigr\}
\end{eqnarray*}
is the event from Lemma~\ref{apriori} and
\begin{eqnarray*}
&& D_j = \bigl\{\bigl|u_{hj}-u_{h(j-1)}-hF(u_{h(j-1)})\bigr|
\leq\mu h\bigr\}
\end{eqnarray*}
is the event from Lemma~\ref{dapriori}, for $\mu>0$ to be chosen. If
$\mu,h>0$ are fixed and $h$ is small enough, then there are constants
$C,\gamma>0$ so that $\mathbb{P}(B_j\cap D_j)\geq1-(C/M)e^{-\gamma
N}$, and since $A_M = \bigcap_{j=1}^M(B_j\cap D_j)$, $\mathbb
{P}(A_M)\geq1-Ce^{-\gamma N}$. For $j=1,\ldots,M$ let
\begin{eqnarray*}
&& E_j = \sup_{\omega\in A_j}\bigl|u_{hj}(\omega) -
u(hj)\bigr|,
\end{eqnarray*}
where $\omega$ denotes an element of the probability space for the
partner model. Letting $u'=F(u)$ denote \eqref{mfeq}, we have
\begin{eqnarray*}
&& u(hj)-u\bigl(h(j-1)\bigr) = \int_{h(j-1)}^{hj}F
\bigl(u(s)\bigr)\,ds.
\end{eqnarray*}
Since $F(u)$ is quadratic in $u$ and its domain is bounded, it is
bounded and Lipschitz continuous, that is, for some $L>0$ and all $u,v$
in the domain, $|F(u)|\leq L$ and $|F(v)-F(u)|\leq L|v-u|$. From the
first inequality, it follows that $|u(s)-u(h(j-1))|\leq L(s-h(j-1))$
for $s \geq h(j-1)$ and from this and the second inequality it follows that
\begin{eqnarray*}
&&\bigl|u(hj)-u\bigl(h(j-1)\bigr)-hF\bigl(u\bigl(h(j-1)\bigr)\bigr)\bigr| \\
&&\qquad =  \biggl
\llvert \int_{h(j-1)}^{hj}\bigl(F\bigl(u(s)\bigr)-F
\bigl(u\bigl(h(j-1)\bigr)\bigr)\bigr)\,ds \biggr\rrvert
\\
&&\qquad \leq \int_{h(j-1)}^{hj}\bigl|F\bigl(u(s)\bigr)-F\bigl(u
\bigl(h(j-1)\bigr)\bigr)\bigr|\,ds
\\
&&\qquad \leq \int_{h(j-1)}^{hj}L\bigl|u(s)-u\bigl(h(j-1)\bigr)\bigr|
\,ds
\\
&&\qquad \leq \int_{h(j-1)}^{hj}L^2(s-hj)\,ds =
L^2\int_0^h s \,ds =
L^2h^2/2.
\end{eqnarray*}
Also,
\begin{eqnarray*}
\bigl|u_{hj}-u(hj)\bigr|
&=& \bigl|u_{hj}-u_{h(j-1)}-hF(u_{h(j-1)})+u_{h(j-1)}-u
\bigl(h(j-1)\bigr)
\\
&&{}+ hF(u_{h(j-1)}) -hF\bigl(u\bigl(h(j-1)\bigr)\bigr) + u\bigl(h(j-1)
\bigr) \\
&&{}+ hF\bigl(u\bigl(h(j-1)\bigr)\bigr)-u(hj)\bigr|
\\
&\leq& \bigl|u_{hj}-u_{h(j-1)}-hF(u_{h(j-1)})\bigr| +
\bigl|u_{h(j-1)}-u\bigl(h(j-1)\bigr)\bigr|
\\
&&{}+ \bigl|hF(u_{h(j-1)})-hF\bigl(u\bigl(h(j-1)\bigr)\bigr)\bigr|
\\
&&{}+ \bigl|u(hj)-u
\bigl(h(j-1)\bigr)-hF\bigl(u\bigl(h(j-1)\bigr)\bigr)\bigr|
\end{eqnarray*}
so using the definition of $A_j$, letting $E_0:=|u_0-u(0)|\leq\delta$
and using once more Lipschitz continuity of $F$ it follows that for $j=1,\ldots,M$,
\begin{eqnarray*}
&& E_j \leq\mu h + E_{j-1} + hLE_{j-1} +
L^2h^2/2 = (1+hL)E_{j-1} + h\bigl(\mu +
hL^2/2\bigr).
\end{eqnarray*}
Setting $q=(1+hL)$ and $r=\mu+hL^2/2$ and iterating the inequality
$E_j\leq qE_{j-1} + hr$, we find\vspace*{1pt} $E_M \leq q^ME_0 + [(q^M-1)/(q-1)]hr
\leq q^M[E_0+hr/(q-1)] = (1+hL)^M[E_0 +hr/(hL)] = (1+LT/M)^M[E_0+r/L]
\leq e^{LT}[E_0+r/L]\leq e^{LT}[\delta+r/L]$ and the same inequality
holds for all $E_j,j=1,\ldots,M$. Since on $A_j$, $|u_s-u_{hj}|\leq C_1h$
for $h(j-1) \leq s \leq hj$, on $A_M$ we find for $j=1,\ldots,M$ and
$h(j-1)\leq s \leq hj$ that
\begin{eqnarray*}
\bigl|u_s-u(s)\bigr| &\leq& |u_s-u_{hj}| +
\bigl|u_{hj}-u(hj)\bigr| + \bigl|u(hj)-u(s)\bigr|
\\
&\leq& C_1h + E_j + Lh \leq h(C_1+L) +
e^{LT}[\delta+r/L]
\end{eqnarray*}
and taking $h,\mu,\delta>0$ small enough, this is $\leq\varepsilon$.
\end{pf}
Our first application of Proposition~\ref{mfest} is to control $y_t$.

\begin{lemma}\label{ydyn}
For each $\varepsilon>0$, there are constants $C,T,\gamma>0$ so that from
any value $y_0\in[0,1]$,
\begin{eqnarray*}
&&\mathbb{P}\Bigl(\sup_{T\leq t \leq e^{\gamma N}}\bigl|y_t-y^*\bigr|\leq
\varepsilon\Bigr)\geq 1-Ce^{-\gamma N}.
\end{eqnarray*}
Moreover, if $|y_0-y^*| \leq2\varepsilon/3$ we may take $T=0$.
\end{lemma}

\begin{pf}
Let $y'=f(y)$ denote the $y'$ equation in \eqref{mfeq} and let $\phi
(t,y)$, $\phi:[0,1]\times\mathbb{R}_+ \rightarrow[0,1]$ denote the
flow for this equation, that is, the unique function satisfying
$\partial_t \phi(t,y) = f(\phi(t,y))$ and $\phi(0,y)=y$ for each
$(t,y)$ in its domain. Since $\phi(t,0)\leq\phi(t,y)\leq\phi(t,1)$ and
$\lim_{t\rightarrow\infty}\phi(t,y)=y^*$ for each $y\in[0,1]$, for
each $\varepsilon>0$ there is $T>0$ so that $|\phi(T,y)-y^*|\leq\varepsilon
/3$ for all $y \in[0,1]$. Letting $y(t)=\phi(t,y_0)$ and using
Proposition~\ref{mfest}, there are constants $C_1,\gamma_1>0$ depending
on $\varepsilon$ but not on $y_0$ so that with probability $\geq
1-C_1e^{-\gamma_1 N}$, $|y_T-y^*| \leq|y_T-y(T)| + |y(T)-y^*| \leq
\varepsilon/3+\varepsilon/3=2\varepsilon/3$. Then, for $t \geq0$ and $y \in
[y^*-(2\varepsilon/3),y^*+(2\varepsilon/3)]$,
\begin{eqnarray*}
&& y^*-(2\varepsilon/3) \leq\phi\bigl(t,y^*-(2\varepsilon/3)\bigr)\leq\phi(t,y)
\leq\phi \bigl(t,y^*+(2\varepsilon/3)\bigr) \leq y^*+(2\varepsilon/3)
\end{eqnarray*}
and since all solutions approach $y^*$ there is $h>0$ so that $\phi
(h,y^*-2\varepsilon/3)\geq y^*-\varepsilon/3$ and $\phi(h,y^*+2\varepsilon/3)
\leq y^*+\varepsilon/3$. Thus, for the given value of $h$ and any solution
$y(t)$ of $y'=f(y)$, if $|y(T)-y^*|\leq2\varepsilon/3$ then
$|y(t)-y^*|\leq2\varepsilon/3$ for $t\geq T$ and $|y(T+h)-y^*|\leq
\varepsilon/3$. Given $y_T$ such that $|y_T-y^*|\leq2\varepsilon/3$ and
setting $y(T)=y_T$, by Proposition~\ref{mfest} there are constants
$C_2,\gamma_2>0$ so that $\sup_{T\leq t \leq T+h}|y_t-y(t)| \leq
\varepsilon/3$ with probability $\geq1-C_2e^{-2\gamma_2 N}$, in which case
\begin{eqnarray*}
\sup_{T \leq t \leq T+h}\bigl|y_t-y^*\bigr|  &\leq & \sup
_{T \leq t \leq
T+h}\bigl|y_t-y(t)\bigr|+\sup_{T \leq t
 \leq   T+h}\bigl|y(t)-y^*\bigr|\\
&\leq & \varepsilon/3 + 2\varepsilon/3 = \varepsilon
\end{eqnarray*}
and $|y_{T+h} - y^*| \leq|y_{T+h} - y(T+h)| + |y(T+h)-y^*| \leq
\varepsilon/3 + \varepsilon/3 = 2\varepsilon/3$ with the same probability.
Iterating this for $e^{\gamma_2 N}/h$ time steps, we find that
\begin{eqnarray*}
&&\sup_{T \leq t \leq e^{\gamma_2 N}}\bigl|y_t - y^*\bigr| \leq\max
_{i \in\{
1,\ldots,e^{\gamma_2 N}/h\}}\sup_{T+(i-1)h \leq t \leq T+ih}\bigl|y_t - y^*\bigr| \leq
\varepsilon
\end{eqnarray*}
with probability $\geq1- (C_2/h)e^{\gamma_2 N}e^{-2\gamma_2 N} = 1 -
(C_2/h)e^{-\gamma_2 N}$, then choose $C=C_1+C_2/h$ and $\gamma= \min
(\gamma_1,\gamma_2)$. Note that if $|y_0-y^*|\leq2\varepsilon/3$, the
iteration step is immediately applicable, in which case we may take $T=0$.
\end{pf}

\section{Macroscopic behaviour}\label{secmacro}
In this section, we prove the macroscopic side of Theorem~\ref{thm1}
that is, when $|V_0|\geq\varepsilon N$. We begin with the analogue of Lemma~\ref{mfmt} for the partner model, which we refer to later on as
monotonicity. As for the MFE, define $ip_t:=si_t+ii_t$.

\begin{lemma}\label{pmmt}
Let $\leq$ denote the partial order on $\mathbb{R}^3$ given by $u\leq
v\Leftrightarrow u_j\leq v_j, \forall j \in\{1,2,3\}$, and let
$u_t=(i_t,ip_t,ii_t)$. If\vspace*{1pt} $(V_t^{(1)},E_t^{(1)})$ and
$(V_t^{(2)},E_t^{(2)})$ are two copies of the partner model with
$E_0^{(1)}=E_0^{(2)}$ and $V_0^{(1)}\subseteq V_0^{(2)}$ then with
respect to the coupling given by the graphical construction,
$E_t^{(1)}=E_t^{(2)}$ and $V_t^{(1)}\subseteq V_t^{(2)}$ for $t>0$ and
correspondingly $y_t^{(1)}=y_t^{(2)}$ and $u_t^{(1)}\leq u_t^{(2)}$.
\end{lemma}

\begin{pf}
If $E_0^{(1)}=E_0^{(2)}$ then $E_t^{(1)}=E_t^{(2)}=:E_t$ for $t>0$.
Given $\{E_t:t\geq0\}$, the only transitions affecting $V_t^{(1)}$ and
$V_t^{(2)}$ are recovery of infectious sites and transmission from
infectious to healthy sites along open edges, both of which preserve
the order $V_t^{(1)}\subseteq V_t^{(2)}$. The equality\vspace*{1pt}
$y_t^{(1)}=y_t^{(2)}$ follows directly from $|E_t^{(1)}|=|E_t^{(2)}|$
and the inequality $u_t^{(1)}\leq u_t^{(2)}$ follows directly from
$V_t^{(1)}\subseteq V_t^{(2)}$.
\end{pf}
Using Proposition~\ref{mfest} and Lemma~\ref{pmmt}, we can prove the
macroscopic part of Theorem~\ref{thm1} when $R_0\leq1$. In this
section, $u_t$ will generally refer to $(i_t,si_t,ii_t)$ or
$(i_t,ip_t,ii_t)$, with $y_t$ written separately.

\begin{proposition}\label{thm1.1}
If $R_0\leq1$, for each $\varepsilon>0$ there are constants $C,T,\gamma
>0$ so that, from any initial configuration, with probability $\geq
1-Ce^{-\gamma N}$, $|V_T|\leq\varepsilon N$.
\end{proposition}

\begin{pf}
By Lemma~\ref{pmmt}, it is enough to show the result holds when $V_0=V$
that is, everyone is initially infectious; in this case $y_0=1-2E_0/N$,
$i_0=y_0$ and $ip_0=ii_0=(1-y_0)/2$. Let $u_t=(i_t,ip_t,ii_t)$ and let
$(y(t),u(t))$ be the solution to the MFE with $y(0)=y_0$ and
$u(0)=u_0$. By Lemma~\ref{ydyn} and Proposition~\ref{mfest}, for each
$\delta>0$ there are constants $C_1,T_1,\gamma_1>0$ so that with
probability $\geq1-C_1e^{-\gamma_1 N}$, $|y_{T_1}-y^*| \leq\delta$
and $|u_{T_1}-u(T_1)|\leq\delta$, so with the same probability
$|(y_{T_1},u_{T_1})-(y^*,u(T_1))| \leq\delta$.

Recall the set $\Lambda^*$ and let $(y^*,\overline{u}(t))$ be the
solution to the MFE with $\overline{u}(0) = (y^*,(1-y^*)/2,(1-y^*)/2)$.
As shown in the proof of Lemma~\ref{mfthmsuff}, $\overline{u}(t)$
decreases to an equilibrium. Since $R_0\leq1$, $(0,0,0)$ is the only
equilibrium, so $\overline{u}(t)\rightarrow(0,0,0)$ as $t\rightarrow
\infty$. Moreover, $\overline{u}(0)\geq v$ for each $v \in\Lambda^*$
so for any solution $(y^*,u(t))$, $\overline{u}(0)\geq u(0)$. By Lemma~\ref{mfmt}, $\overline{u}(t)\geq u(t)$ for $t\geq0$, so there is $T_2$
not depending on $u(0)$ so that $|u(T_2)| \leq\varepsilon/2$. Using
Proposition~\ref{mfest}, there are constants $C_2,\gamma_2,\delta>0$
not depending on $u(0)$ so that with probability $\geq1-C_2e^{-\gamma
_2 N}$, if $|(y_0,u_0)-(y^*,u(0))|\leq\delta$ then $|u_{T_2}| \leq
|u(T_2)|+|u_{T_2}-u(T_2)| \leq\varepsilon/2+\varepsilon/2=\varepsilon$. Letting
$T=T_1+T_2$, $C=C_1+C_2$ and $\gamma=\min(\gamma_1,\gamma_2)$ and
combining the two steps completes the proof.
\end{pf}
Using similar ideas, we can prove the macroscopic part of Theorem~\ref
{thm1} when $R_0>1$. Before showing the approach to equilibrium, we
first have to show long time survival of the infection, and to do that
we need the following result concerning the MFEs.

\begin{lemma}\label{liftup}
Suppose $R_0>1$ and let $v \in\mathbb{R}^3$ with $|v|=1$ be an unstable
eigenvector of the MFEs on $\Lambda^*$ as given in the proof of Theorem~\ref{mfthm}, written in $(i,ip,ii)$ coordinates. For $0<\delta' \leq
\delta$ let $(y(t),u(t))$ be a solution to the MFE with $|y(0)-y^*|\leq
\delta$ and $u(0):=(i(0),ip(0),ii(0)) = \delta' v$. If $\delta>0$ is
small enough, then there is $T>0$ so that $\min_j u_j(T)\geq2\delta'$
for all $0< \delta' \leq\delta$.
\end{lemma}

\begin{pf}
First, write the MFE \eqref{mfeq2}, without the $y$ equation, in matrix
form as follows:
%
\begin{equation}
\label{mfmtx} %
\pmatrix{i'
\cr
ip'
\cr
ii'} = %
\pmatrix{ -(1+r_+y) & r_- & r_-
\cr
r_+
\bigl(y^*-i/2\bigr) & -(1+r_-) & 1
\cr
r_+i/2 & \lambda& -(2+r_-+\lambda)}
\pmatrix{ i
\cr
ip
\cr
ii}.
\end{equation}
The $y$ dynamics\vspace*{1pt} proceeds as in \eqref{mfeq}, and note $|y(t)-y^*|\leq
|y(0)-y^*|$ for $t>0$. Write \eqref{mfmtx} as $u' = A(i,y)u$ with
$u=(i,si,ii)^{\top}$ to emphasize the dependence on $i,y$. As noted in
the proof of Theorem~\ref{mfthm}, if $R_0>1$ then $A:= A(0,y^*)$ has a
positive eigenvalue $\mu>0$ with positive eigenvector $v$ such that
$|v|=1$, so the system $v' = Av$ has solutions $v(t)=cve^{\mu t}$ for
any $c>0$. Let $|\cdot|$ denote the operator norm and let
\begin{eqnarray*}
&& L=\sup_{(i,y)\in[0,1]^2}\bigl|A(i,y)\bigr|
\end{eqnarray*}
then any solution $u(t)$ to \eqref{mfmtx} has $|u(t)|\leq|u(0)|e^{Lt}$
for $t>0$. Fix $T>0$, then for each $\varepsilon>0$, by continuity there
is $\delta>0$ so that if $\max( y-y^*,i)\leq e^{LT}\delta$ then
$|A(i,y)-A|\leq\varepsilon$. Let $|y(0)-y^*|\leq\delta$ and for $0<\delta
'\leq\delta$ let $u(t)$ be the solution to~\eqref{mfmtx} with
$u(0)=\delta' v$, then for $0\leq t \leq T$,
\begin{eqnarray*}
\bigl|(u-v)'\bigr| &=& \bigl|A(i,y)u - Av\bigr| \leq\bigl|\bigl(A(i,y)-A\bigr)u\bigr| +
\bigl|A(u-v)\bigr|
\\
&\leq& \bigl|A(i,y)-A\bigr||u| + |A||u-v|
\\
&\leq& \varepsilon|u|+ L|u-v|
\\
&\leq& \varepsilon e^{Lt}\delta' + L|u-v|.
\end{eqnarray*}
Letting $v(0)=u(0)$, defining $E(t) := |u(t)-v(t)|$, noting that
$E(0)=0$ and integrating,
\begin{eqnarray*}
&& E(T) \leq e^{LT}\varepsilon\delta' T.
\end{eqnarray*}
Since $v(T) = \delta v e^{\mu T}$,
\begin{eqnarray*}
\min_j u_j(T) &\geq& \min
_j v_j(T) - \max_j
\bigl|v_j(T)-u_j(T)\bigr|
\\
&\geq& \delta' e^{\mu T}\min_j
v_j - \varepsilon e^{LT}\delta'T
\\
&=& \delta' e^{\mu T}\Bigl(\min_j
v_j -\varepsilon e^{(L-\mu)T}T\Bigr).
\end{eqnarray*}
Choose $T>0$ so that $e^{\mu T}\min_j v_j/2 \geq2$, then choose
$\varepsilon>0$ so that $\varepsilon e^{(L-\mu)T}T \leq\min_j v_j/2$, then
it follows that $\min_j u_j(T) \geq2\delta'$.
\end{pf}
Now we can show long-time survival of the infection when $R_0>1$ and
$|V_0|\geq\varepsilon N$.

\begin{lemma}\label{stayup}
Suppose $R_0>1$. For each $\varepsilon>0$, there are constants $\delta
,C,\gamma>0$ so that if $|V_0|\geq\varepsilon N$ then
\begin{eqnarray*}
&& \mathbb{P}\Bigl(\inf_{0\leq t \leq e^{\gamma N}}|V_t|\geq\delta N
\Bigr)\geq 1-Ce^{-\gamma N}.
\end{eqnarray*}
\end{lemma}

\begin{pf}
Recall that an event holds with high probability or w.h.p. in $N$ if for
$N$ large enough it occurs with probability $\geq1-Ce^{-\gamma N}$ for
some $C,\gamma>0$. If $|V_0|\geq\varepsilon N$ then $\max
(i_0,ip_0,ii_0)\geq\varepsilon/3$, so in view of Lemma~\ref{pmmt} it is
enough to prove the result starting from $u_0:=(i_0,ip_0,ii_0) \in
\mathcal{E} := \{(\varepsilon/3,0,0),\break (0,\varepsilon/3,0),(0,0,\varepsilon/3)\}$.
For $\delta_1>0$, by Proposition~\ref{ydyn} there are $T,\gamma_1>0$ so
that w.h.p. $|y_t-y^*|\leq\delta_1$ for $T\leq t \leq e^{\gamma_1 N}$. If
$u(0)\neq(0,0,0)$ then for $t>0$, $\min_j u_j(t)>0$; this is shown for
$u(0)\in\Lambda^*$ in the proof of Lemma~\ref{mfthmsuff}, but the same
proof applies if $y\neq y^*$. Also, since $(0,0,0)$ is an equilibrium
solution, by uniqueness of solutions $u(t)\neq(0,0,0)$ for $0\leq t
\leq T$, so by continuity of solutions $\inf\{|u(t)|:0\leq t \leq T\}
>0$. Therefore, there exists $0<\delta_2\leq\delta_1$ so that $\min_j
u_j(T) \geq\delta_2$ and $\inf\{\max_j u_j(t):0 \leq t \leq T\} \geq
\delta_2$ for all $u(0)\in\mathcal{E}$. For $u_0=u(0)\in\mathcal{E}$
with $y_0=y(0) \in[0,1]$, by Proposition~\ref{mfest}, w.h.p. $|u_t-u(t)|
\leq\delta_2/2$ for $0\leq t \leq T$ in which case $\min
(i_T,ip_T,ii_T)\geq\delta_2/2$ and $\inf\{\max(i_t,ip_t,ii_t):0\leq t
\leq T\}\geq\delta_2/2$, which means that for the eigenvector $v$ with
$|v|=1$ mentioned in the proof of Lemma~\ref{liftup}, $(i_T,si_T,ii_T)
\geq(\delta_2/2)v$, and also $|V_t|\geq(\delta_2/2)N$ for $0\leq t
\leq T$.

Taking $y(t)=y_t$ and $u(T) = (\delta_2/2)v$, if $|y_t-y^*|\leq\delta
_1$ then by Lemma~\ref{liftup} there is $h>0$ so that $\min_j u_j(T+h)
\geq \delta_2$, and as before there is $\delta_3>0$ so that $\inf\{\max_j u_j(t):T \leq t \leq T+h\}\geq\delta_3$. By Lemma~\ref{pmmt} and the
last paragraph, it is enough to consider the case $u_T=u(T)=(\delta
_2/2)v$. Letting $\delta=\min(\delta_2/2,\delta_3/2)$ and using
Proposition~\ref{mfest}, with probability $\geq1-Ce^{\gamma_2 N}$,
$|u_{t}-u(t)|\leq\delta$ for $T\leq t \leq T+h$, in which case
$u_{T+h}\geq(\delta_2/2)v$ and $|V_t|\geq N\min(i_t,ip_t,ii_t)\geq
(\delta_3/2)N$ for $T\leq t \leq T+h$. Letting $\gamma=\min(\gamma
_1/2,\gamma_2/2)$ and iterating for $e^{\gamma N}/h$ time steps as in
the proof of Lemma~\ref{ydyn}, w.h.p. $|V_t| \geq N \min
(i_t,ip_t,ii_t)\geq(\delta_3/2)N$ for $T \leq t \leq e^{\gamma N}$.
Combining with the previous estimate, w.h.p. $|V_t|\geq\delta N$ for
$0\leq t \leq e^{\gamma N}$ as we wanted to show.
\end{pf}
We now wrap up the macroscopic side of Theorem~\ref{thm1}.

\begin{proposition}\label{thm1.2}
Suppose $R_0>1$ and let $(y^*,i^*,ip^*,ii^*)$ with $i^*>0$ be the
nontrivial equilibrium solution to the MFE \eqref{mfeq2}. Let
$u_t=(i_t,ip_t,ii_t)$ and let $u^*=(i^*,ip^*,ii^*)$. For each $\varepsilon
>0$, there are constants $C,T,\gamma>0$ so that if $|V_0|\geq\varepsilon
N$ then
\begin{eqnarray*}
&& \mathbb{P}\Bigl(\sup_{T \leq t \leq e^{\gamma N}}\bigl|(y_t,u_t)-
\bigl(y^*,u^*\bigr)\bigr|\leq \varepsilon\Bigr)\geq1-Ce^{-\gamma N}.
\end{eqnarray*}
\end{proposition}

\begin{pf}
We begin with the lower bound. As shown in the proof of Lemma~\ref
{stayup} there are $T_1,h_1,\delta_1,\gamma_1>0$ so that w.h.p. $\min
(i_t,ip_t,ii_t)\geq\delta_1$, and thus $u_t\geq\delta_1 v$, for
$t=T_1+kh_1$, $k=1,\ldots,(e^{\gamma_1 N}-T_1)/h_1$, where $v$ with
$|v|=1$ is the eigenvector from Lemma~\ref{liftup}. Let $y(0)=y^*$ and
$u(0):=(i(0),ip(0),ii(0))=\delta_1 v$. If $\delta_1>0$ is small enough,
then $u_j'(0)>0$ in each coordinate and since $u^*\neq(0,0,0)$ is
unique, as shown in the proof of Lemma~\ref{mfthmsuff} $u(t)$ is
increasing with respect to $(i,ip,ii)$ coordinates and $\lim_{t\rightarrow\infty}u(t)=u^*$, and in particular $u(t)\leq u^*$ for
$t\geq0$. We will need the stronger fact $u_j(t)<u_j^*$ for $j=1,2,3$.
Looking to the equations for $i',ip',ii'$ in \eqref{mfeq2}, the
derivative of each variable increases with the other two variables, and
of course is equal to $0$ at $u^*$. If we had $i(t)=i^*$, then since
$ip(t)\leq ip^*$ and $ii(t)\leq ii^*$ we would have $i'<0$ which
contradicts the fact that $u(t)$ is increasing, and the same applies to
$ip(t)$ and $ii(t)$.

Using the above facts, there is $T_2$ so that $u(T_2)\geq u^*-\varepsilon
/2$, and since $0<\min_j (u_j^*-u_j(T_2))=:\varepsilon'\leq\varepsilon$,
there is $h_2$ so that $u(T_2+h_2)\geq u^*-\varepsilon'/2$. By Proposition~\ref{mfest}, there is $\delta_2>0$ so that if $u_0=u(0)$ and
$|y_0-y^*|\leq\delta_2$ then w.h.p. $|u_t-u(t)|\leq\varepsilon'/2$ for
$T_2\leq t \leq T_2+h_2$ in which case $u_t\geq u^*-\varepsilon$ for
$T_2\leq t \leq T_2+h_2$ and $u_{T_2+h_2}\geq u^*-\varepsilon'$, which
means\vspace*{1pt} that $u_{T_2+h_2}\geq u(T_2)$. By Lemma~\ref{ydyn}, there are
$T_3,\gamma_2$ so that w.h.p. $|y_t-y^*|\leq\min(\delta_2,\varepsilon)$ for
$T_3\leq t \leq e^{\gamma_2 N}$. Let $k$ be such that $T_1+kh_1 \geq
T_3$ and let $T_4=T_1+kh_1$, then setting $u(T_4)=\delta_1 v$, w.h.p.
$u_{T_4} \geq u(T_4)$ so it is enough to consider the case where
$u_{T_4}=u(T_4)$. Letting $T=T_4+T_2$, then for some $\gamma_3>0$, with
probability $\geq1-Ce^{-\gamma_3 N}$, $u_t\geq u^*-\varepsilon$ for
$T\leq t \leq T+h_2$ and $u_{T+h_2}\geq u(T)$. Letting $\gamma=\min
(\gamma_2/2,\gamma_3/2)$ and iterating for $(e^{\gamma N}-T)/h_2$ time
steps (subtracting $T$ to make sure $y_t$ stays in bounds) as in the
proof of Lemma~\ref{ydyn} it follows that $u_t\geq u^*-\varepsilon$ for
$T\leq t \leq e^{\gamma N}$.

To prove the upper bound, it is enough to consider any value of $y_0$
and let $u_0 = (y_0,(1/2)(1-y_0),(1/2)(1-y_0))$. Setting $y(0)=y^*$ and
$u(0)=(y^*,(1/2)(1-y^*),(1/2)(1-y^*))$, then as shown in the proof of
Lemma~\ref{mfthmsuff}, $u(t)$ decreases to~$u^*$. Moreover,
$u_j(t)-u_j^*>0$ for the same reason as above, so there is $T_1$ so
that $u(T_1)\leq u^*+\varepsilon/2$, and since $0<\min_j(u_j(T_1)-u_j^*)=:\varepsilon'\leq\varepsilon$, there is $h$ so that
$u(T_1+h)\geq u^*-\varepsilon'/2$. By Proposition~\ref{mfest}, there is
$\delta>0$ so that if $\max(|u_0-u(0)|,|y_0-y(0)|)\leq\delta$ then
w.h.p.
$|u_t-u(t)|\leq\varepsilon'/2$ for $T_1\leq t \leq T_1+h$ in which case
$u_t \leq u^*+\varepsilon$ for $T_1\leq t \leq T_1+h$ and $u_{T_1+h}\leq
u^*+\varepsilon'$ which means that $u_{T_1+h}\leq u(T_1)$. By Lemma~\ref
{ydyn}, there are $T_2,\gamma_1$ so that w.h.p. $|y_t-y^*|\leq\delta$ for
$T_2\leq t \leq e^{\gamma_1 N}$. Letting $T=T_1+T_2$ and setting
$u(T_2)=(y^*,(1/2)(1-y^*),(1/2)(1-y^*)$ and
$u_{T_2}=(y_{T_2},(1/2)(1-y_{T_2}),(1/2)(1-y_{T_2}))$, then for some
$\gamma_2>0$, with probability $\geq1-Ce^{-\gamma_2 N}$, $u_t \leq
u^*+\varepsilon$ for $T\leq t \leq T+h$ and $u_{T+h} \leq u(T)$. Letting
$\gamma= \min(\gamma_1/2,\gamma_2/2)$ and iterating\vspace*{1pt} for $(e^{\gamma
N}-T)/h$ time steps it follows as before that $u_t\leq u^*+\varepsilon$
for $T \leq t \leq e^{\gamma N}$.
\end{pf}
In the next section, we use a comparison to prove that if $R_0<1$ the
infection disappears quickly from the population. To make this work, we
will need a complementary result to Lemma~\ref{stayup}.

\begin{lemma}\label{staydown}
If $R_0\leq1$, then for each $\varepsilon>0$ there are $C,T,\gamma>0$ so that
\begin{eqnarray*}
&& \mathbb{P}\Bigl(\sup_{T\leq t \leq e^{\gamma N}}|V_t|\leq\varepsilon
N\Bigr) \geq 1-Ce^{-\gamma N}.
\end{eqnarray*}
\end{lemma}

\begin{pf}
The proof is similar to that of Lemma~\ref{stayup}. Letting $\overline
{u}=(y^*,(1-y^*)/2,(1-y^*)/2)$ as in Lemma~\ref{mfthmsuff} and letting
$(y^*,\overline{u}(t))$ be the solution to the MFE with $\overline
{u}(0)=\overline{u}$, since $\overline{u}(t)$ decreases to $(0,0,0)$
and $\overline{u}\geq v$ for all $v \in\Lambda^*$, there is $T_1$ so
that for any solution $(y^*,u(t))$, $|u(T_1)|\leq\varepsilon/6$, and
since $\varepsilon':=\min_ju_j(T_1)>0$, there is $h$ so that
$|u(T_1+h)|\leq\varepsilon'/2$. There is $\delta>0$ so that if $\max
(|u_0-u(0)|,|y_0-y^*|)\leq\delta$ then w.h.p. $|u_t-u(t)|\leq\min(\varepsilon
'/2,\varepsilon/6)$ for $T_1\leq t \leq T_1+h$ in which case $|u_t|\leq
\varepsilon/3$ for $T_1\leq t \leq T_1+h$ and $|u_{T_1+h}|\leq\varepsilon'$
which means $u_{T_1+h}\leq u(T_1)$. There are $\gamma_1,T_2>0$ so that
w.h.p. $|y_t-y^*|\leq\delta$ for $T_2\leq t \leq e^{\gamma_1 N}$. By
monotonicity, it is enough to consider $u_{T_2}=\overline{u}$. Letting
$u(T_2)=u_{T_2}$ and $T=T_1+T_2$, there are $C_1,\gamma_2$ so that with
probability $\geq1- C_1e^{\gamma_2 N}$, $|u_t|\leq\varepsilon/3$ for
$T\leq t \leq T+h$ and $u_T\leq u(T)$. Letting $\gamma=\min(\gamma
_1,\gamma_2)$ and iterating for $(e^{\gamma N}-T)/h$ time steps, w.h.p.
$|u_t|\leq\varepsilon/3$ and thus $|V_T|\leq\varepsilon N$ for $T \leq t
\leq e^{\gamma N}$.
\end{pf}

\section{Microscopic behaviour}\label{secmicro}
In this section, we compare the partner model in the regime $|V| \leq
\varepsilon N$ for small $\varepsilon>0$ to a branching process to get
decisive information when $R_0\neq1$.

\subsection{Subcritical case: $R_0<1$}
First, we introduce the comparison process to use when $R_0<1$.

\begin{definition}\label{ubpdef}
Define the \emph{upperbound process} (UBP) $B_t = (\mathcal
{I}_t,\mathcal{SI}_t,\mathcal{II}_t)$ on state space $\{0,1,2,\ldots\}^3$
with parameter $0 \leq\delta\leq y^*$ by the following transitions:
\begin{itemize}
\item$\mathcal{I}\rightarrow\mathcal{I}-1$ at rate $\mathcal{I}$,
\item$\mathcal{I}\rightarrow\mathcal{I}-1$ and $\mathcal
{SI}\rightarrow\mathcal{SI}+1$ at rate $r_+(y^*-\delta)\mathcal{I}$,
\item$\mathcal{SI}\rightarrow\mathcal{SI}+1$ at rate $2r_+\delta
\mathcal{I}$,
\item$\mathcal{II}\rightarrow\mathcal{II}+1$ at rate $r_+\delta
\mathcal{I}$,
\item$\mathcal{SI}\rightarrow\mathcal{SI}-1$ at rate $\mathcal{SI}$,
\item$\mathcal{SI}\rightarrow\mathcal{SI}-1$ and $\mathcal
{I}\rightarrow\mathcal{I}+1$ at rate $r_- \mathcal{SI}$,
\item$\mathcal{SI}\rightarrow\mathcal{SI}-1$ and $\mathcal
{II}\rightarrow\mathcal{II}+1$ at rate $\lambda\mathcal{SI}$,
\item$\mathcal{II}\rightarrow\mathcal{II}-1$ and $\mathcal
{SI}\rightarrow\mathcal{SI}+1$ at rate $2 \mathcal{II}$,
\item$\mathcal{II}\rightarrow\mathcal{II}-1$ and $\mathcal
{I}\rightarrow\mathcal{I}+2$ at rate $r_- \mathcal{II}$.
\end{itemize}
\end{definition}

Note the UBP describes the evolution of the total number of particles
of each of the three types $\mathcal{I},\mathcal{SI},\mathcal{II}$ in a
multi-type continuous-time branching process; for an introduction to
branching processes, see \cite{bp}. We now show that for fixed $R_0<1$,
if $\delta>0$ is small enough the UBP quickly dies out.

\begin{lemma}\label{ubpdown}
For fixed $\lambda,r_+,r_-$, let $B_t$ denote the UBP with parameter $\delta'$
and let $R_0$ be as defined in \eqref{r0eq0}. If $R_0<1$,
there are $C,\delta>0$ so that if $|B_0|\leq N$ and $\delta'\leq\delta
$ then
\begin{eqnarray*}
&& \mathbb{P}\bigl(|B_{C\log N}|=0\bigr) \rightarrow1 \qquad\mbox{as }N\rightarrow
\infty.
\end{eqnarray*}
\end{lemma}

\begin{pf}
For a multi-type continuous time branching process $B_t =
(b_1(t),\break \ldots, b_n(t))$, with $b_j(t)$ denoting the number of type $j$
particles alive at time $t$, we can extract some useful information
from the \emph{mean matrix} $M_t$ defined by $m_{ij}(t) = \mathbb
{E}(b_j(t) \vert b_k(0)=\delta_{ik})$. Since particles evolve
independently, $\mathbb{E}(B_t) = B_0M_t$ and it is not hard to show
that $M_t$ satisfies the equation
\begin{eqnarray*}
&&\frac{d}{dt}M_t = AM_t
\end{eqnarray*}
and, therefore, $M_t = \exp(At)$, where $A = (r_{ij})$ is the matrix
whose entries $r_{ij}$ give the rate at which a particle of type $i$
produces particles of type $j$. If $\real(\lambda)<0$ for each
eigenvalue $\lambda$ of $A$, then letting $\gamma_0 = \min\{|\real
(\lambda)|:\lambda\in\sigma(A)\}$ where $\sigma(\cdot)$ denotes the
spectrum, from standard matrix theory it follows that for any $\gamma
_1<\gamma_0$, there is $C_1>0$ so that $m_{ij}\leq C_1e^{-\gamma_1t}$
for each pair $ij$. Since each $b_i(t)$ is valued on nonnegative integers,
\begin{eqnarray*}
\mathbb{P}\bigl(B_t\neq(0,\ldots,0)\bigr)&\leq&\sum
_i \mathbb{P}\bigl(b_i(t) \neq0\bigr) \leq
\sum_i\mathbb{E}b_i(t)
\\
&= &\sum_{ij}b_i(0)m_{ij}(t)
\leq\bigl|B(0)\bigr|n^2C_1e^{-\gamma_1 t}.
\end{eqnarray*}
If $|B(0)|\leq N$, then\vspace*{1pt} letting $t=C\log N$ for $C>1/\gamma_1$ and
setting $\gamma= C\gamma_1-1$ and $C_2=n^2C_1$ we find
\begin{eqnarray*}
 \mathbb{P}\bigl(B_{C\log N}\neq(0,\ldots,0)\bigr)&\leq &
NC_2e^{-\gamma_1 C\log N} = NC_2N^{-\gamma_1 C}\\
& =&
C_2N^{1-\gamma_1 C} = C_2N^{-\gamma}
\end{eqnarray*}
which tends to $0$ as $N\rightarrow\infty$. In our case,
\begin{eqnarray*}
&& A = A(\delta) = %
\pmatrix{-\bigl(1+r_+\bigl(y^*-\delta\bigr)\bigr)
& r_+\bigl(y^*+\delta\bigr) & r_+\delta
\cr
r_- & -(1+r_-+\lambda) & \lambda
\cr
2r_- & 2 & -(2+r_-) }.
\end{eqnarray*}
Letting $\sigma(A)$ denote the spectrum and defining the \emph{spectral
abcissa} $\mu(A):=\max\{\real(\lambda):\lambda\in\sigma(A)\}$, if $\mu
(A(\delta))<0$, then the real part of each eigenvalue of $A$ is
negative, and the above argument applies. By continuity of eigenvalues
in the entries of a matrix, it is enough to show $\mu(A(0))<0$, since
then there is $\delta>0$ so that if $\delta'\leq\delta$ then $\mu
(A(\delta')) \leq\mu(A(0))/2 < 0$. Setting $\delta=0$,
\begin{eqnarray*}
&& A(0)= %
\pmatrix{-\bigl(1+r_+y^*\bigr) & r_+y^* & 0
\cr
r_- &
-(1+r_-+\lambda) & \lambda
\cr
2r_- & 2 & -(2+r_-)}
\end{eqnarray*}
and looking to Section~\ref{secmf} we see that $A(0,0)$ is the
(transpose of the) linearized matrix at $(0,0,0)$ for the MFE on
$\Lambda^*$, which we denote $A$. As noted in Remark~\ref{r0remark},
$(0,0,0)$ is locally asymptotically stable when $R_0<1$, and in the
proof of Theorem~2 in \cite{watm} this is done by showing that $\mu(A)<0$.
\end{pf}

We now complete the proof of the case $R_0<1$ in Theorem~\ref{thm1}.

\begin{proposition}\label{microdown}
If $R_0<1$ there are constants $C,T,\gamma>0$ so that, from any initial
configuration,
\begin{eqnarray*}
&& \mathbb{P}\bigl(|V_{T+C\log N}| =0\bigr)\rightarrow1 \qquad\mbox{as } N\rightarrow
\infty.
\end{eqnarray*}
\end{proposition}

\begin{pf}
Let $U_t := (I_t,\mathit{SI}_t,\mathit{II}_t)$ denote variables in the partner model and
for $\delta>0$ such that $y^*-\delta\geq0$ and $y^*+\delta\leq1$,
let $B_t$ denote the UBP with parameter~$\delta$. We first describe a
coupling with the property that $U_0\leq B_0\Rightarrow U_t\leq B_t$
for $t>0$, with respect to the usual partial order $U \leq V
\Leftrightarrow U_j\leq V_j,j=1,2,3$. For $j=1,\ldots,10$, define a
countable number of independent Poisson point processes (p.p.p.'s) $\{
e_j(n):n=1,2,\ldots\}$ with respective rates $1,r_+,r_+,1,r_-,\lambda
,2,r_-,r_+,r_-$, together with independent uniform $[0,1]$ random
variables attached to each event in $e_2(n),e_3(n),e_9(n),
n=1,2,\ldots.$
These correspond to the nine transitions listed in the definition of
the UBP, except that the second and third transition in the UBP are
lumped into $e_2$, plus an additional transition for $S+S\rightarrow
\mathit{SS}$ and one for $\mathit{SS}\rightarrow S+S$. Note that the rates of
$e_2,e_3,e_9$ appear too large at the moment and are corrected in the
next paragraph.

Construct the UBP one transition at a time as follows, letting
$(\mathcal{I},\mathcal{SI},\mathcal{II})$ denote the present state.
Each event in $e_1(1),\ldots,e_1(\mathcal{I})$ reduces $\mathcal{I}$ by
$1$. For an event in $e_2,e_3$ let $p$ denote the corresponding uniform
$[0,1]$ random variable. If an event in $e_2(1),\ldots,e_2(\mathcal{I})$
occurs and $p \leq(y^*-\delta)$, reduce $\mathcal{I}$ by 1 and
increase $\mathcal{SI}$ by 1, while if $y^*-\delta< p \leq y^*+\delta$
simply increase $\mathcal{SI}$ by~1. If an event in
$e_3(1),\ldots,e_3(\mathcal{I})$ occurs and $p\leq\delta$, increase
$\mathcal{II}$ by 1. Each event in $e_4(1),\ldots,e_4(\mathcal{SI})$
reduces $\mathcal{SI}$ by 1, each event in $e_5(1),\ldots,e_5(\mathcal
{SI})$ reduces $\mathcal{SI}$ by 1 and increases $\mathcal{I}$ by 1,
each event in $e_6(1),\ldots,e_6(\mathcal{SI})$ reduces $\mathcal{SI}$ by
1 and increases $\mathcal{II}$ by 1, each event in
$e_7(1),\ldots,e_7(\mathcal{II})$ reduces $\mathcal{II}$ by 1 and
increases $\mathcal{SI}$ by 1, and each event in $e_8(1),\ldots,e_8(\mathcal
{II})$ reduces $\mathcal{II}$ by 1 and increases $\mathcal{I}$ by~2. It
can be checked that the transition rates are correct.

Similarly, construct the Markov chain $(S_t,I_t,\mathit{SS}_t,\mathit{SI}_t,\mathit{II}_t)$ for
the partner model as follows, letting $(S,I,\mathit{SS},\mathit{SI},\mathit{II})$ denote the
present state. Define $\alpha_t = y_t-y^*-i_t$ and $\beta_t =
i_t/2-1/(2N)$ and note that $\alpha_t$ and $\beta_t$ are piecewise
constant in time. Each event in $e_1(1),\ldots,e_1(I)$ reduces $I$ by $1$
and increases $S$ by 1. If an event in $e_2(1),\ldots,e_2(I)$ occurs and
$p\leq y^*+\alpha_t$ reduce $S$ and $I$ by 1 and increase $\mathit{SI}$ by~1. If
an event in $e_3(1),\ldots,e_3(I)$ occurs and $p \leq\beta_t$ reduce $I$
by 2 and increase $\mathit{II}$ by 1. Each event in $e_4(1),\ldots,e_4(\mathit{SI})$ reduces
$\mathit{SI}$ by 1 and increases $\mathit{SS}$ by 1, each event in $e_5(1),\ldots,e_5(\mathit{SI})$
reduces $\mathit{SI}$ by 1 and increases $S$ and $I$ by 1, and events in
$e_6,e_7,e_8$ have the same effect as before. If an event in
$e_9(1),\ldots,e_9(S)$ occurs and $p \leq s_t/2-1/(2N)$ reduce $S$ by 2
and increase $\mathit{SS}$ by 1, and each event in $e_{10}(1),\ldots,e_{10}(\mathit{SS})$
reduces $\mathit{SS}$ by 1 and increases $S$ by 2. Recalling that
$U_t:=(I_t,\mathit{SI}_t,\mathit{II}_t)$, if $U_0 \leq B_0$ and $\sup_{s \leq t}\max
(|\alpha_s|,\beta_s)\leq\delta$ then $U_t \leq B_t$ since (as can be
easily checked) the order is preserved at each transition.

By Lemma~\ref{ydyn}, there are $T_1,\gamma_1>0$ so that w.h.p.
$|y_t-y^*|\leq\delta/2$ for $T_1\leq t \leq e^{\gamma N}$ and since
$R_0<1$, by Lemma~\ref{staydown} there are $T_2,\gamma_2$ so that
$|V_t|\leq(\delta/2)N$, and thus $i_t \leq\delta/2$ for $T_2\leq t
\leq e^{\gamma_2 N}$. Letting $T=\max(T_1,T_2)$ and $\gamma=\min(\gamma
_1,\gamma_2)$, w.h.p. $\max(|\alpha_t|,\beta_t)\leq\delta$ for $T\leq t
\leq e^{\gamma N}$. Setting $B_T=U_T$ and using Lemma~\ref{ubpdown}
completes the proof.
\end{pf}

\subsection{Supercritical case: $R_0>1$}

We introduce the comparison process for $R_0>1$, which is similar to
the UBP, but different.

\begin{definition}\label{lbpdef}
Define the \emph{lowerbound process} (LBP) $B_t = (\mathcal
{I}_t,\mathcal{SI}_t,\mathcal{II}_t)$ on state space $\{0,1,2,\ldots\}^3$
with parameters $\delta\geq0$ such that $y^*-\delta\geq0$ by the
following transitions:
\begin{itemize}
\item$\mathcal{I}\rightarrow\mathcal{I}-1$ at rate $(1+2r_+\delta
)\mathcal{I}$,
\item$\mathcal{I}\rightarrow\mathcal{I}-1$ and $\mathcal
{SI}\rightarrow\mathcal{SI}+1$ at rate $r_+(y^*-\delta)\mathcal{I}$,
\item$\mathcal{I}\rightarrow\mathcal{I}-2$ at rate $r_+\delta\mathcal{I}$,
\item$\mathcal{SI}\rightarrow\mathcal{SI}-1$ at rate $\mathcal{SI}$,
\item$\mathcal{SI}\rightarrow\mathcal{SI}-1$ and $\mathcal
{I}\rightarrow\mathcal{I}+1$ at rate $r_- \mathcal{SI}$,
\item$\mathcal{SI}\rightarrow\mathcal{SI}-1$ and $\mathcal
{II}\rightarrow\mathcal{II}+1$ at rate $\lambda\mathcal{SI}$,
\item$\mathcal{II}\rightarrow\mathcal{II}-1$ and $\mathcal
{SI}\rightarrow\mathcal{SI}+1$ at rate $2 \mathcal{II}$,
\item$\mathcal{II}\rightarrow\mathcal{II}-1$ and $\mathcal
{I}\rightarrow\mathcal{I}+2$ at rate $r_- \mathcal{II}$.
\end{itemize}
\end{definition}

As before, the LBP describes the evolution of the total number of
particles of each of the three types $\mathcal{I},\mathcal{SI},\mathcal
{II}$ in a multi-type continuous-time branching process. We now show
that for fixed $R_0>1$, if $\delta>0$ is small enough then the LBP survives.

\begin{lemma}\label{lbpup}
Let $B_t$ denote the $LBP$ with parameter $\delta'$. If $\lambda
,r_+,r_-$ are such that $R_0>1$ then there are $C,\delta>0$ so that if
$\delta'\leq\delta$ then\break $\liminf_{N\rightarrow\infty}\mathbb{P}(B_{C
\log N}\neq(0,0,0))>0$ and
\begin{eqnarray*}
&& \mathbb{P}\bigl(|B_{C \log N}|\geq\delta N \vert B_{C \log N}
\neq(0,0,0)\bigr) \rightarrow1 \qquad\mbox{as }N \rightarrow\infty.
\end{eqnarray*}
\end{lemma}

\begin{pf}
As in the proof of Lemma~\ref{ubpdown}, define the mean matrix
$M(t)=\exp(At)$ and the spectral abcissa $\mu(A)$. If $\delta'=0$ for
both the UBP and the LBP they coincide, in which case $A$ is the
transpose of the linearized matrix at $(0,0,0)$ of the MFE on $\Lambda
^*$. As shown in the proof of Theorem~\ref{mfthm}, if $R_0>1$ then $\mu
(A)>0$. By continuity of eigenvalues in the entries of a matrix, there
is $\delta>0$ so that if $\delta'\leq\delta$ then $\mu(A(\delta'))\geq
\mu(A)/2 >0$. As shown in V.7 of \cite{bp}, if $M(t)$ is such that for
some $t_0>0$ and each entry $m_{ij}(t)$ of $M(t)$ one has
$m_{ij}(t_0)>0$ (which is the case here), then $\mu(A)=:\lambda_1$ is
an eigenvalue of $A$, and if $\lambda_1>0$ the process is said to be
\emph{supercritical}. In this case, $B_te^{-\lambda_1 t} \rightarrow
Wv$ where $v$ is a left eigenvector of $A$ with eigenvalue $\lambda_1$
and $W$ is a real-valued random variable. Setting $t=C\log N$ with
$C>1/\lambda_1$ and letting $\gamma= C\lambda_1>1$, $B_{C \log
N}N^{-\gamma} \rightarrow Wv$, so for each $\varepsilon>0$,
\begin{eqnarray*}
&& \liminf_{N\rightarrow\infty}\mathbb{P}\bigl(|B_{C\log N}|\geq\delta N\bigr)
\geq \lim_{N\rightarrow\infty}\mathbb{P}\bigl(|B_{C\log N}| \geq
\varepsilon N^{\gamma}\bigr) = \mathbb{P}\bigl(W|v|\geq\varepsilon\bigr)
\end{eqnarray*}
and letting $\varepsilon\rightarrow0^+$ and using continuity of measure,
\begin{eqnarray*}
&& \liminf_{N\rightarrow\infty}\mathbb{P}\bigl(|B_{C\log N}|\geq\delta N\bigr)
\geq \mathbb{P}(W>0).
\end{eqnarray*}
Under a mild regularity assumption on the offspring distribution that
holds trivially in this case, $\mathbb{P}(W>0)=\lim_{t\rightarrow\infty
}\mathbb{P}(B_t\neq(0,0,0))>0$. Since $|B_t| \geq\delta N$ implies
$B_t\neq(0,0,0)$, this means $\limsup_{N\rightarrow\infty}\mathbb
{P}(|B_{C\log N}|\geq\delta N) \leq\break \lim_{t\rightarrow\infty}\mathbb
{P}(B_t\neq(0,0,0)) =\mathbb{P}(W>0)$, so $\lim_{N\rightarrow\infty
}\mathbb{P}(|B_{C\log N}|\geq\delta N)$ exists and is equal to $\mathbb
{P}(W>0)$. The result then follows by observing that for $t,x>0$,
$\mathbb{P}(|B_t| \geq x \vert B_t \neq(0,0,0)) = \mathbb{P}(|B_t|\geq
x)/\mathbb{P}(B_t \neq(0,0,0))$.
\end{pf}
We now complete the proof of Theorem~\ref{thm1}.

\begin{proposition}
If $R_0>1$, there are constants $\delta,p,C,T>0$ so that if $|V_0|>0$
then $\mathbb{P}(|V_{T+C\log N}| \geq\delta N )\geq p$.
\end{proposition}

\begin{pf}
We use the same approach as in the proof of Proposition~\ref
{microdown}. Let $U_t := (I_t,\mathit{SI}_t,\mathit{II}_t)$ denote variables in the
partner model and for $\delta_1>0$ such that $\delta_1\leq1$,
$y^*-\delta_1\geq0$ and $y^*+\delta_1\leq1$, let $B_t$ denote the LBP
with parameter $\delta_1$. Let $e_1,\ldots,e_{10}$ be as in the proof of
Proposition~\ref{microdown}.

Construct the LBP one transition at a time as follows, letting
$(\mathcal{I},\mathcal{SI},\mathcal{II})$ denote the present state.
Each event in $e_1(1),\ldots,e_1(\mathcal{I})$ reduces $\mathcal{I}$ by
$1$. For an event in $e_2,e_3$ let $p$ denote the corresponding uniform
$[0,1]$ random variable. If an event in $e_2(1),\ldots,e_2(\mathcal{I})$
occurs and $p \leq(y^*-\delta_1)$, reduce $\mathcal{I}$ by 1 and
increase $\mathcal{SI}$ by 1, while if $y^*-\delta_1 < p \leq y^*+\delta
_1$ simply reduce $\mathcal{I}$ by 1. If an event in
$e_3(1),\ldots,e_3(\mathcal{I})$ occurs and $p\leq\delta_1$, reduce
$\mathcal{I}$ by 2. Events in $e_4,e_5,e_6,e_7,e_8$ have the same
effect as in the dynamics of the UBP. The Markov chain
$(S_t,I_t,\mathit{SS}_t,\mathit{SI}_t,\mathit{II}_t)$ for the partner model is constructed in the
same way as in the proof of Proposition~\ref{microdown}, with $\alpha
_t,\beta_t$ defined in the same way, and it is easy to check in this
case that if $U_0\geq B_0$ and $\sup_{s \leq t}\max(|\alpha_s|,\beta
_s)\leq\delta_1$ then $U_t \geq B_t$.

Define the stopping time $\tau= \inf\{t:|U_t|\geq\delta_1 N/2\}$ and
note that $|V_{\tau}|\geq(\delta_1/2) N$. By Lemma~\ref{stayup} and
using the strong\vspace*{1pt} Markov property, there are $\delta,\gamma>0$ so that
w.h.p. $|V_t|\geq\delta N$ for $\tau\leq t \leq\tau+e^{\gamma N}$. There
are $T,\gamma>0$ so that w.h.p. $|y_t - y^*|\leq\delta_1/2$ for $T\leq t
\leq e^{\gamma N}$. If $\tau\leq T$, then since $T$ is fixed, we are
done. If $t<\tau$ then $i_t \leq\delta_1/2$, so letting $B_T=U_T$, if
$T\leq t <\tau$ then $\max(|\alpha_t|,\beta_t)\leq\delta_1$, so
$U_t\geq B_t$ for $T \leq t < \tau$. The result follows from this
inequality and from Lemma~\ref{lbpup}.
\end{pf}

\section*{Acknowledgements}
The authors wish to thank Chris Hoffman for the suggestion to study the
model on the complete graph, as well as the referee for a thorough
reading of the article and for helpful comments.

%





\printaddresses
\end{document}